\magnification 1200
\input amssym.def
\input amssym.tex
\parindent = 40 pt
\parskip = 12 pt
\font \heading = cmbx10 at 12 true pt
 at 22 true pt
\font \medheading =cmbx7 at 14 true pt
 at 7 true pt
\def \R{{\bf R}}

\centerline{\medheading $L^p$ boundedness of maximal averages over hypersurfaces in $\R^3$}
\rm
\line{}
\line{}
\line{}
\centerline{\heading Michael Greenblatt}
\line{}
\centerline{August 23, 2010}
\baselineskip = 12 pt
\font \heading = cmbx10 at 14 true pt
\line{}
\line{}
\noindent{\bf 1. Introduction and statement of results}

\vfootnote{}{This research was supported in part by NSF grant DMS-0919713}Let $Q$ be a smooth hypersurface 
in $\R^3$ and let $q_0$ be a point on $Q$. Let $d\sigma(q)$ denote the standard surface measure on $Q$.
For a small cutoff function $\phi$ supported near $q_0$, we define the maximal operator $M$, initially on 
Schwarz functions, by 
$$Mf(x) = \sup_{t > 0} |\int_Q f(x - tq) \phi(q)\,d\sigma(q)| \eqno (1.1)$$
Our goal is to determine for which $p$ is the maximal operator $M$ bounded on $L^p$. Note that by subadditivity
of maximal operators of the form $(1.1)$, one can prove $L^p$ boundedness of the nonlocalized analogue of 
$(1.1)$ over a compact surface by doing a partition of unity and reducing to $(1.1)$. 
If $L$ is an invertible linear transformation and $M_L$ denotes the maximal operator corresponding to the
surface $L(Q)$, then one can easily check from the definitions that $M_L f (x) = |det(L)| M(f \circ L)
(L^{-1}x)$. Hence when studying $L^p$ boundedness properties one may always replace $M$ by $M_L$. In particular, we may assume that
near $q_0$, $Q$ is given as the graph of a function $g(x,y)$ such that if $(x_0,y_0)$ denotes the point
in the $x$-$y$ plane which $q_0$ lies above, then $\nabla g(x_0,y_0) = 0$. 

The earliest work in this area was done in the case where $Q$ is an $n$-dimensional sphere in $\R^{n+1}$,
when Stein [St1] showed $M$ is bounded on $L^p$ 
iff $p > {n + 1 \over n}$ for $n > 1$. This was later generalized by Greenleaf [Gr] to surfaces of 
nonvanishing Gaussian curvature, with some further results for when the Hessian has rank between 1 and $n$.
The $n = 1$ case was later proven by Bourgain [B]. In [So] Sogge showed in any dimension that whenever $Q$
has at least one nonvanishing principal curvature, $M$ is bounded on $L^p$ for all $p > 2$. The case of 
convex surfaces of finite line type has been extensively analyzed; we refer to [IoSe2] and [NaSeWa] for more
information on these situations.

Although there are
many interesting issues when $p \leq 2$, for the purposes of this paper we always assume $p > 2$. $M$ is 
trivially bounded on $L^{\infty}$, and 
if $M$ is bounded on some $L^p$, by interpolating with the $L^{\infty}$ case one has that $M$ is bounded on
$L^{p'}$ for $p' > p$. Hence our goal is to determine the optimal $p_0 \geq 2$ for which $M$ is bounded on 
$L^p$ for $p > p_0$. If $Q$ is tangent to the tangent plane $T_{q_0}(Q)$ to infinite 
order at $q_0$, then as long as $0 \notin T_{q_0}(Q)$ a relatively straightforward argument shows that
$M$ is unbounded on $L^p$ for all finite $p$. Conversely, if the Gaussian curvature of $Q$ does not vanish
to infinite order at $q_0$, then by [SoSt] $M$ is bounded on $L^p$ for some finite $p$. (See [CoMa] for
another result of this kind.) It is entirely possible that in any dimension, $M$ is bounded on some 
finite $L^p$ whenever $Q$ is not tangent to $T_{q_0}(Q)$ to infinite order at $q_0$. Hence in general 
we expect $p_0$ to be finite. 

\noindent {\bf Definition 1.1.} Let $d(x,y)$ denote the vertical distance between $Q$
and $T_{q_0}(Q)$ above $(x,y)$. The {\it height} $h(q_0)$ is defined 
to be the reciprocal of the supremum of all $\epsilon$ for which 
the integral of $|d(x,y)|^{-\epsilon}$ is finite on at least one neighborhood of $q_0$.

For the case $n = 2$ considered in this paper, a good $L^p$ boundedness theorem for $p > 2$ was proven in 
[IkKeMu2]. Their main theorem can be stated as follows.

\noindent {\bf Theorem [IkKeMu2]}. Suppose the origin is not contained in $T_{q_0}(Q)$. If $\phi(q)$
is supported on a sufficiently small neighborhood of $q_0$ then $M$ is 
bounded on $L^p$ for $p > \max(h(q_0), 2)$. When $h(q_0) \geq 2$ and $\phi(q_0) \neq 0$ this exponent is 
sharp in that $M$ is unbounded on $L^p$ for $p < h(q_0)$ and if $Q$ is real-analytic,
then $M$ is unbounded on $L^{h(q_0)}$ as well. 

The purpose of this paper is to provide a relatively short alternative approach to the $p > 2,$ $n = 2$ situation by
extending the methods of [G1] and using facts about the adapted coordinate systems described below.
There will once again be exceptional cases not covered, but due 
to the differences in the methods the exceptional cases will be quite different from and not mutually 
exclusive to the exceptional cases of [IkKeMu2], which occur when $0 \in T_{q_0}(Q)$.

\noindent {\bf Newton Polygons and Adapted Coordinates.}

\noindent We now give some relevant terminology which will be used throughout this paper. Below, $R(x,y)$
denotes a smooth function defined on a neighborhood of the origin with nonvanishing Taylor expansion at
the origin.

\noindent {\bf Definition 1.2.} Let $R(x,y) = \sum_{a,b} R_{ab}x^ay^b$ denote the 
Taylor expansion of $R(x,y)$ at the origin. For any $(a,b)$ for which $R_{ab} \neq 0$, let $Q_{ab}$ be the
quadrant $\{(x,y) \in \R^2: 
x \geq a, y \geq b \}$. Then the {\it Newton polygon} $N(R)$ of $R(x,y)$ is defined to be 
the convex hull of the union of all $Q_{ab}$.  

In general, a Newton polygon consists of finitely many (possibly zero) bounded edges of negative slope
as well as an unbounded vertical ray and an unbounded horizontal ray.

\noindent {\bf Definition 1.3.} The {\it Newton distance} $d(R)$ of $R(x,y)$ is defined to be 
$\inf \{t: (t,t) \in N(R)\}$.

One often uses $(t_1,t_2)$ coordinates to write equations of lines relating
to Newton polygons, so as to distinguish from the $x$-$y$ variables of the domain of $R(x,y)$. The line 
in the $t_1$-$t_2$ plane with equation $t_1 = t_2$ comes up so frequently it has its own name:

\noindent {\bf Definition 1.4.} The {\it bisectrix} is the line in the $t_1$-$t_2$ plane with equation
$t_1 = t_2$.

A key role in the above theorems as well as our theorems to follow is played by the following polynomials.

\noindent {\bf Definition 1.5}. Suppose $e$ is a compact edge of $N(R)$. Define $R_e(x,y)$
by $R_e(x,y) = \sum_{(a,b) \in e} R_{ab} x^{a}y^{b}$. In other words $R_e(x,y)$ is the sum of the terms
of the Taylor expansion of $f$ corresponding to $(a,b) \in e$. 

\noindent {\bf Definition 1.6}. Suppose $R(x,y)$ has nonvanishing Taylor expansion at the origin such that 
$R(0,0) = 0$ and $\nabla R(0,0) = 0$. Then $R(x,y)$ is said to be in {\it nonadapted coordinates} if the 
bisectrix intersects $N(R)$ in the interior of a compact edge $e$ of $N(R)$ such that $R_e(1,y)$ has a zero 
of order greater than $d(R)$. If $R(x,y)$ is not in nonadapted coordinates, then $R(x,y)$ is said to be 
{\it adapted coordinates}.

The significance of adapted and nonadapted coordinates was first discovered by Varchenko [V] for the 
real-analytic case and for the general smooth case by Ikromov-M\"uller [IkMu]. Namely, define $\epsilon_0$ to be 
the supremum of all $\epsilon$ for which $|R|^{-\epsilon}$ is integrable in at least one neighborhood of
the origin. Equivalently, $\epsilon_0$ is the supremum of the epsilon such that on some neighborhood of 
$(0,0)$ one has $|\{(x,y): |R(x,y)| < t \}| < C t^{\epsilon}$ for some $C$.
Then by [V] and [IkMu] one has $d(R) \leq {1 \over \epsilon_0}$, with equality holding if and only if 
$R(x,y)$ is in
adapted coordinates. Furthermore, one has the following. Suppose $R(x,y)$ is not in adapted coordinates and
let $e$ be the edge of $N(R)$ intersecting the bisectrix in its interior. By [V] and [IkMu], if the slope 
$m_e$ of $e$ is at least $-1$, then there is a smooth $\psi(x)$ with $\psi(0) = 0$
such that $R(x,y + \psi(x))$ is in adapted coordinates. By switching the roles of the $x$ and $y$ axes,
this means that if $m_e \leq -1$, there is a smooth $\psi(y)$ with $\psi(0) = 0$ such that $R(x + \psi(y)
,y)$ is in adapted coordinates. Thus [V] and [IkMu] show that there necessarily is a "nice" 
coordinate change after which the growth rate of $R$ at the origin is given in the above way by 
${1 \over d(R)}$ and where $R$ is in adapted coordinates; these facts in turn are used in [IkKeMu2] in their 
proof of the $L^p$ boundedness properties of $M$ in their theorem above. Another useful aspect of adapted
coordinates proven in [V] and [IkMu] is that if one is in adapted coordinates, then the order of any zero
of any $R_e(1,y)$ (cf. Definition 1.5) is at most $d(R)$. 

Let $a$ denote the order of the zero of $R(x,y)$ at $(0,0)$. Then for a generic linear transformation $T$, 
$\partial_y^{a}(R \circ T)(0,0)$ and $\partial_x^{a}(R \circ T)(0,0)$ are both
nonzero. Thus the Newton polygon of $R \circ T$ is entirely on or above the line $t_1 + t_2 = a$. But
no point of $N(R \circ T)$ can be below this line; otherwise $r \circ T$ would have a zero of order less than
$a$ at the origin. We conclude that $N(R \circ T)$ is exactly $\{(t_1,t_2) \in \R^2: t_1 + t_2 \geq
a\}$. Since the compact edge of $N(R \circ T)$ has slope $-1$, by the above discussion there is a 
$\psi(x)$ such that $R \circ T(x, y + \psi(x))$ is in adapted coordinates. Note that the Newton polygon of
$R \circ T(x, y + \psi(x))$ still has its upper vertex at $(0, a)$ and that the slope of each edge 
of this Newton polygon is at least -1. This motivates the following definition.

\noindent {\bf Definition 1.7.} Suppose $R(0,0)= 0$ and $\nabla R(0,0) = 0$. Then $R(x,y)$ is said to be in
{\it generic adapted coordinates} if $R(x,y)$ is in adapted coordinates, each edge of $N(R)$ has slope at
least $-1$, and $N(R)$ intersects the $y$-axis at some point $(0,a)$. 

Note that the above definition implies that $R(x,y)$ has a zero of order $a$ at the origin. We now are 
in a position to state the main theorem of this paper. Recall we are working in the situation where $Q$
is a surface with a distinguished point $q_0 = (x_0,y_0,z_0)$, such that near $q_0$ the surface $Q$ is the
graph of some $g(x,y)$ with $g(x_0,y_0) = z_0$ and $\nabla g(x_0,y_0) = (0,0)$. Let $f(x,y) = g(x + x_0,
y + y_0) - z_0$, and define $D(x,y)$ to be the Hessian determinant of $f$ at $(x,y)$. Then our main
theorem is as follows

\noindent {\bf Theorem 1.1.} Suppose $M$ is as defined in $(1.1)$. If $\phi(q)$ is supported on a 
sufficiently small neighborhood of $q_0$, then $M$ is bounded on $L^p$ for $p >
\max(h(q_0), 2)$ as long as neither of the following two exceptional situations occurs.

\noindent {\bf a)} $D(x,y)$ has a zero of infinite order at $(0,0)$.

\noindent {\bf b)} Whenever $T$ is an invertible linear transformation and $\psi(x)$ a smooth function with 
a zero of order $b > 0$ at the origin such that $F(x,y) = (f \circ T) (x, y + \psi(x))$ is in generic 
adapted coordinates, then the bisectrix intersects $N(F)$ in the interior of a compact edge $e$ with slope
$m_e$ with $|m_e| <  {1 \over b}$ such that $p(y) = \partial_y (F_e(1,y))$ or $\partial_y (F_e(-1,y))$ has a zero 
of order greater than $\max(1,d(F)-1)$ at some $y_0$ for which $p(y_0) \neq 0$.

Here $F_e(x,y)$ is as in Definition 1.5. If $\psi(x)$ has a zero of infinite order at $x = 0$, then we take 
${1 \over b} = 0$ in Theorem 1.1b); in other words, the exceptional condition of b) will not be satisfied. 
As is explained in [IkKeMu2], using a theorem in [IoSa1] one can show the exponent $h(q_0)$ is best possible
when $q_0 \notin T_{q_0}(Q)$, assuming $\phi(q_0) \neq 0$. If $q_0 \in T_{q_0}(Q)$, then sometimes one can 
do better as the maximal operator starts to resemble a traditional Hardy-Littlewood maximal operator in two 
dimensions. 

The exceptional
situation a) of this paper is necessitated by our use of damping functions in conjunction with the theorem 
of Sogge and Stein [SoSt] that we will describe below. It is unlikely that it can be avoided without using
substantial additional ideas. A canonical example of when the first kind of exception situation occurs is
when $F(x,y) = p(ax + by)$, in other words, when $F(x,y)$ is effectively a function of one variable.

Exceptional situation b) may be viewed as rare in the sense that it requires a 
certain polynomial to have a zero of high order among other things; however, a simple example of where it happens is the
function $(y + x^a)^b + x^c$ for $c \geq  b > 2$, $a > 1$, and $ab < c$. Although we will not prove it here,
further manipulations involving Newton polygons can be used to show that if the exceptional condition holds
for one of $\partial_y (F_e(1,y))$ and $\partial_y (F_e(-1,y))$, it holds for the other. Hence in the statement
of Theorem 1.1 we could have just used one of the two functions.

The exceptional situation b) arises for the following reason. The main theorem
of [G1] gives Theorem 1.1 if $f(x,y)$ is already in adapted coordinates. Much of the analysis of [G1] 
carries over even in nonadapted coordinates; the $y$-variable shift by $\psi(x)$ does not interfere with
most of the argument. The exception occurs when one cannot avoid using integrations by parts in the
$x$-variable in arguments resembling the proof of the Van der Corput lemma. This happens when the bisectrix
intersects $N(F)$ in the interior of a bounded edge satisfying
the above conditions on $F_e(1,y)$ and its $y$ derivatives. If the order $b$ of the zero of $\psi(x)$ at the origin is at least
${1 \over |m_e|}$, then $\psi(x)$ is too small to cause any serious problems in such integrations by parts.
If $b$ is less than ${1 \over |m_e|}$ then it seems difficult to adapt the 
arguments of [G1] to these situations. So as long as we can find some linear $T$ such that $f \circ T$
avoids such situations, then the methods of [G1] can be adapted to the current situation and Theorem 1.1
can be proved. 

The stipulation that $b < {1 \over |m_e|}$ is more than a technical improvement in the statement of Theorem
1.1. For example, if $f(x,y)$ is a mixed homogeneous function, then it is not hard to show
that in converting to adapted coordinates $\psi(x)$ can always be taken to be of the form $cx^b$, and 
if the bisectrix of $F(x,y) = f(x, y + cx^b)$ intersects $N(F)$ in the interior of a compact edge of slope
$m_e$ then $b = {1 \over |m_e|}$. Hence Theorem 1.1 covers the mixed homogeneous case (as long as $D(x,y)$
does not have a zero of infinite order at the origin) because we include the condition. For $h(q_0) \geq 2$
this result was first proved in [IkKeMu1].

\noindent {\bf Strategy of the Proof.}

Ever since [SoSt] one successful method of proving $L^p$ boundedness theorems for maximal operators such
as $(1.1)$ has involved embedding $M$ in an analytic family $M_z$. The idea is that one replaces the standard
surface measure $d\sigma(q)$ with the damped surface measure $e^{z^2}|h|^z d\sigma(q)$, and then defines the maximal
operator $M_z$ to be the analogue of $(1.1)$ with $d\sigma(q)$ replaced by $e^{z^2}|h|^z d\sigma(q)$. One
shows that for some $s_0 < 0$, $M_z$ is bounded on $L^{\infty}$ whenever $Re(z) > s_0$, uniformly in 
$|Im(z)|$ for fixed $Re(z)$, and that for some $s_1 > 0$, $M_z$ is bounded on $L^2$ whenever $Re(z) > s_1$, 
again uniformly in $|Im(z)|$ for fixed $Re(z)$. The $e^{z^2}$ factor is present to ensure uniform $L^2$ 
bounds on each vertical line. Then by a well known interpolation technique for maximal operators
(see Ch 11 of [St2] for details), one obtains a $p_0 > 2$ such that $M = M_0$ is bounded on $L^p$ for $p > p_0$. 
The hope is that the damping function $h(z)$ can be chosen so that $p_0$ is optimal. 

In the above interpolation, the $L^{\infty}$ bounds are typically obtained using the observation that
$||M_z||_{L^{\infty} \rightarrow L^{\infty}}$ is bounded by the $L^1$ norm of $C|h|^{Re(z)}$. Thus as long
as $|h|^{s_0}d\sigma(q)$ is a finite measure, one obtains the desired uniform $L^{\infty}$ bounds on any
vertical line $Re(z) = s$ for $s > s_0$.
For the $L^2$ boundedness, we will use the following consequence of a theorem of Sogge and Stein:

\noindent {\bf Theorem 1.2. [SoSt].} Suppose the surface measure $\phi(q)d\sigma(q)$ is as in $(1.1)$
and there are $C, \epsilon > 0$ such that for all multiindices $\alpha$ with $|\alpha| = 0,1$ the Fourier
transform of the measure $d\sigma_z(q) =  e^{z^2}|h(x)|^z \phi(q)d\sigma(q)$ satisfies  
$$|\partial^{\alpha} \hat{\sigma_z} (\lambda)| < C(1 + |\lambda|)^{-{1 \over 2} - \epsilon} \eqno (1.2)$$
Then there is a constant $C'$ depending on $C$ and $\epsilon$ such that $||M_zf||_2 \leq C'||f||_2$ 
for all $f \in L^2$. 

In practice, if one has $(1.2)$ for $\alpha = 0$, it will generally automatically hold for all the first 
derivatives since the effect of taking such a derivative is to replace the cutoff function by another one.
In the proof of Theorem 1.1 in this paper, we will work in this framework. We will first write
the surface $Q$ near $q_0$ as the union of
finitely many "slivers" containing $q_0$ on their boundaries. These slivers will be defined using
the Newton polygon of the function $f(x,y)$ defined above Theorem 1.1 when put in generic adapted 
coordinates; we will effectively be doing a coarse resolution of singularities using Lemmas 2.2 and 2.3 of
the next section. We will then define the damping function separately on each sliver. The damping functions 
will be analogues of those of [G1]. There will be
four types of slivers and showing $(1.2)$ holds will be done by separately estimating the contribution
to $\hat{\sigma_z}(\lambda)$ coming from each type of sliver, again using methods analogous to those of 
[G1]. It should be pointed
out that the idea of dividing a neighborhood of a point into slivers with respect to a Newton polygon
for the purpose of proving oscillatory integral estimates such as $(1.2)$ is quite old; it appears in 
[PSt] and its predecessors for example and was also used in analyzing such maximal operators in 
[IkKeMu1] and [IkKeMu2].

\noindent {\bf 2. Lemmas about Newton polygons; subdivisions into slivers}

In the analysis of this paper, for the appropriate $F(x,y)$ one treats on similar footing all compact edges 
of $N(F)$ intersecting the set $\{(t_1,t_2): t_2 > t_1\}$. To avoid exceptional situations such as those 
of part b) in the statement of Theorem 1.1 for any such edge not intersecting the bisectrix in its interior,
we have the following lemma.

\noindent {\bf Lemma 2.1.} Suppose $R(x,y)$ is a smooth function on a neighborhood of the origin that is
in generic adapted coordinates. Let $d^*$ denote $\max(2, d(R))$. Suppose $e$ is any compact edge of $N(R)$ 
lying entirely on or above the bisectrix. Then if there is $y_0 \neq 0$ such that if $R_e(1,y_0) \neq 0$ and 
$\partial_y R_e(1,y)$ has a zero at $y_0$ of order greater than $d^* - 1$, then $e$ has slope $-1$ and upper
vertex lying on the $y$ axis.

\noindent {\bf Proof.} Since $R(x,y)$ is in generic adapted coordinates, the uppermost vertex of $N(R)$ is
$(0,a)$ for some $a \geq 2$. The point $(d(R),d(R))$ is on $N(R)$, and the line of slope $-1$ containing 
this point intersects $N(R)$ at $(0,2d(R))$. Since all edges of $N(R)$ have slope at least -1 and 
$(d(R),d(R)) \in N(R)$, we conclude that $a \leq 2d(R)$ with equality holding iff there is an edge of
$N(R)$ connecting $(d(R),d(R))$ to $(0,2d(R))$. Hence if $e$ is a compact edge of $N(R)$ lying entirely
on or above the bisectrix, either $e$ is the segment $(d(R),d(R))$ to $(0,2d(R))$ or $a < 2d(R)$. Thus
if $e$ satisfies the assumptions of this lemma and we can show that $R_e(1,y)$ has degree at least $2d(R)$,
then in particular $e$ has slope $-1$ and upper vertex lying on the $y$ axis as needed.

So assume $e$ is a compact edge lying entirely on or above the bisectrix such that there is $y_0 \neq 0$ 
with $R_e(1,y_0) \neq 0$ and such that $\partial_y R_e(1,y)$ has a zero at $y_0$ of order greater than 
$d^*$.  First consider the case where $e$'s lower vertex is $(d(R),d(R))$. In particular, $d(R)$ is an 
integer. We will show that $R_e(1,y)$ has degree at least $2d(R)$. Note that $R_e(1,y)$ has a zero at 
$y = 0$ of order $d(R)$ and therefore 
$\partial_y R_e(1,y)$ has a zero at $y = 0$ of order $d(R) - 1$. Since we are assuming it also has a zero
at $y_0$ of order greater than $d(R) - 1$, $\partial_y R_e(1,y)$ must have degree
least $2d(R) - 1$. Hence $R_e(1,y)$ has degree at least $2d(R)$. This means that the upper vertex
of $e$ is at least $2d(R)$ and by the previous paragraph we are done.

Next, we consider the case where $e$ does not contain $(d(R),d(R))$. Hence $e$ lies entirely above the 
bisectrix, and $R_e(1,y)$ can be written as $y^{d'}p(y)$ where $d' > d(R)$. This means that $\partial_y
R_e(1,y)$ is of the form $y^{d'-1}q(y)$. Since $\partial_y R_e(1,y)$ is assumed to have a zero $y_0 \neq 0$
of order greater than $d(R) - 1$, we can write 
$$\partial_y R_e(1,y) = y^{d'-1}(y - y_0)^{d''-1}r(y) \eqno (2.1)$$
Here $d'' > d(R)$. Note that $(2.1)$ implies that the degree of $\partial_y R_e(1,y)$ is at least $d' + d''
- 2 > 2d(R) - 2$, so the degree of $R_e(1,y)$ is greater than $2d(R) - 1$. Hence the upper vertex of $e$ must be
the upper vertex of $N(R)$; otherwise $N(R)$ would have a vertex at height greater than $2d(R)$ and since
$a \leq 2d(R)$ this can't happen. We conclude the upper vertex of $e$ is given by $(0,a)$ for some $2d(R) 
- 1 < a \leq 2d(R)$. As a result, the degree of $\partial_y R_e(1,y)$ is greater than $2d(R) - 2$ and
at most $2d(R) - 1$.

If $r$ had positive degree, then the degree of $\partial_y R_e(1,y)$ would be at least $d' + d'' - 1 >
2d(R) - 1$, contradicting the above. So $r(y)$ is constant and we may write
$$\partial_y R_e(1,y) = cy^{d'-1}(y - y_0)^{d''-1} \eqno (2.2)$$
Similarly, if $d'$ or $d''$ were equal to $d(R) + 1$ or greater then the degree of $\partial_y R_e(1,y)$ would be
greater than $2d(R) - 1$, again giving a contradiction. So we have $d(R) < d', d'' < d(R) + 1$. Since $d'$ and 
$d''$ are both integers, this means $d' = d''$. "Homogenizing" $(2.2)$, for some $k$ we get that 
$$\partial_y R_e(x,y) = cy^{d'-1}(y - y_0x^k)^{d'-1} \eqno (2.3)$$
Looking at the term of $(2.3)$ whose degree in $y$ is second-highest we see that $k$ is an integer. 
Hence $e$ is an edge of
slope ${-{1 \over k}}$ whose upper vertex is $(0, 2d' - 1)$ and our objective is to show that $k = 1$. Assume 
$k \geq 2$; we will arrive at a contradiction. Since we are dealing with the case that $e$ lies entirely 
above the bisectrix, the point $(d(R),d(R))$ lies above the line containing $e$. Since this line intersects 
the bisectrix at $({k \over k +1}(2d' - 1), {k \over k +1}(2d' - 1))$, we conclude that
$${k \over k +1}(2d' - 1) < d(R) \eqno (2.4)$$
Since $d' > d(R)$ and ${k \over k + 1} \geq {2 \over 3}$, this in turn implies that
$${2 \over 3}(2d(R) - 1) < d(R) \eqno (2.5)$$
Equivalently, $d(R) < 2$. Since $d' < d(R) + 1$ and $d' - 1$ is an integer at least one, $d' - 1 = 1$ and 
$(2.3)$ just becomes
$$\partial_y R_e(x,y) = cy(y - y_0x^k) \eqno (2.6)$$
Now note that $\partial_y R_e(1,y)$ no longer has a zero of order greater than one, so it no longer falls 
under the assumptions of this lemma. So this situation cannot happen; we have arrived at a contradiction
and we are done.

\noindent {\bf Lemma 2.2.} Suppose $R(x,y)$ is a smooth function on a neighborhood of the origin such that
$R(0,0) = 0$. Suppose $v$ is a vertex of $N(R)$ that is the intersection of
compact edges $e_1$ and $e_2$ with slopes $0 > m_1 > m_2$. Let $M_1 = -{1 \over m_1}$ and $M_2 = 
-{1 \over m_2}$.  Let $R_{cd}x^cy^d$ denote 
the term of the Taylor expansion of $R(x,y)$ at the origin corresponding to $v$. Then on a sufficiently small
neighborhood of the origin, there is an $N > 0$ such that if $N |x|^{M_1} < |y| < {1 \over N}|x|^{M_2}$,
then we have 
$${1 \over 2}|R_{cd} x^c y^d| < |R(x,y)| <  2|R_{cd} x^c y^d|\eqno (2.7)$$
If $v$ lies on the $y$ axis and is the upper vertex of a compact edge $e_1$ with slope $-{1 \over M_1}$, then
there is an $N > 0$ such that if $N |x|^{M_1} < |y|$ then once again $(2.7)$ holds.

\noindent {\bf Proof.} Without loss of generality, we may restrict our attention to $(x,y)$ in the upper
right quadrant. Write the Taylor expansion of $R(x,y)$ at the origin as $\sum_{a,b} R_{ab}x^ay^b$. We first 
prove $(2.7)$ in the case where $v$ is the intersection of two compact edges, whose equations we denote by
$t_1 + M_1 t_2 = \alpha_1$ and $t_1 + M_2 t_2 = \alpha_2$. For a large $K$ we can write
$$R(x,y) -  R_{cd}x^cy^d = \sum_{(a,b): c \leq a < M,\,\,d \leq b < M,\,\,(a,b) \neq (c,d)}R_{ab}x^ay^b$$
$$ + \sum_{(a,b):a < c,\,\,d < b < M ,\,\,a + M_2b \geq \alpha_2}R_{ab}x^ay^b + 
\sum_{(a,b): c < a < M,\,\,b < d,\,\,a + M_1b \geq \alpha_1}R_{ab}x^ay^b + E_K(x,y) \eqno (2.8)$$
Here $E_K(x,y)$ satisfies 
$$|E_K(x,y)| < C(|x|^K + |y|^K) \eqno (2.9)$$
 We start by noting 
that the first sum in $(2.8)$ is less than ${1 \over 8}|r_{cd}|x^cy^d$ in absolute value if $(x,y)$ is in a 
sufficiently small neighborhood of the origin, which we may assume. As for the second sum, if one changes coordinates from $(x,y)$ to 
$(x,y')$, where $y' = x^{M_2}y$, then $(x,y') \in [0,1] \times [0,{1 \over N}]$ whenever $y < {1 \over N}
x^{M_2}$. Observe that under this coordinate change, a given term $R_{ab}x^ay^b$ of the second sum becomes
$R_{ab}x^{a + M_2b}(y')^b$. Since
$a + M_2b \geq \alpha_2$ and $b > d$ in each term in the second sum, the entire sum can be written as
$x^{\alpha_2}(y')^d (y'f(x,y'))$ for some $f(x,y')$ which is a polynomial in $y'$ and a fractional power of
$x$. Thus if $N$ is sufficiently large, whenever  $y' < {1 \over N}$ the sum is of absolute value 
less than ${1 \over 8}|R_{cd}|x^{\alpha}(y')^d = {1 \over 8}|R_{cd}|x^cy^d$. Since $y' < {1 \over N}$ is
equivalent to $y < {1 \over N}x^{M_2}$, these are the bounds we need.

The third sum is dealt with in exactly the same way, reversing the roles of the $x$ and $y$ axes and the
edges $e_1$ and $e_2$. Lastly,
since ${1 \over N}x^{M_2} > y > Nx^{M_1}$ the error term $E_K(x,y)$ is less than ${1 \over 8} |r_{cd}|x^cy^d$ in 
absolute value for small $|x|, |y|$ if $K$ is chosen sufficiently large. Putting these all together,
we get that $|R(x,y) - R_{cd}x^cy^d| < {1 \over 2}|R_{cd}|x^cy^d$ as needed. This completes the proof of
Lemma 2.2 for the case where $(c,d)$ is the intersection of two compact edges of $N(R)$.

We now move to the case where $(c,d)$ is on the $y$-axis and is the upper vertex of a compact edge $e_1$ of
$N(R)$. We again examine the sum $(2.8)$. In the case at hand, since $(c,d)$ is on the $y$-axis, the second
sum in $(2.8)$ is empty. The second sum is where the condition $y < {1 \over N}x^{M_2}$ was used above,
and the third sum is where the condition $y > Nx^{M_1}$ was used. Since we have no second sum, the lack of
a condition $y < {1 \over N}x^{M_2}$ holding does not cause any problem in repeating the above argument. For
the third sum we use the condition $y > Nx^{M_1}$ exactly as before, and for the first and fourth sum the 
previous argument works unmodified. Hence $(2.7)$ holds again and we are done.

\noindent {\bf Lemma 2.3.} Suppose $R(x,y)$ is a smooth function of $y$ and a fractional power of $x$ on a neighborhood of the origin such that
$R(0,0) = 0$. Write the Taylor expansion of $R(x,y)$ at the origin as $\sum_{a,b}
R_{ab}x^ay^b$. For a given $M > 0$ let $R_M(x,y)$ denote the sum of the nonzero terms of this Taylor
expansion for which $a + Mb$ is minimized; in particular $R_M(x,y)$ is either of the form $R_e(x,y)$ for
a compact edge $e$ of $N(R)$ or is equal to $R_{cd}x^cy^d$ for a vertex $(c,d)$ of $N(R)$. Denote this 
minimal value of $a + Mb$ by $\alpha$. Then for any
$r \in \R$ and any $\epsilon > 0$, there is a $\delta > 0$ such that on the set $\{(x,y) \in \R^2: 0 < x < 
\delta, (r - \delta)x^M < y < (r + \delta)x^M\}$ we have
$$|R(x,y) - R_M(x,y)| < \epsilon x^{\alpha} \eqno (2.10)$$
\noindent {\bf Proof.} On the region $\{(x,y) \in \R^2: 
0 < x < \delta , (r - \delta)x^M < y < (r + \delta)x^M\}$, we do the coordinate change $(x,y) = (x,x^My')$, 
converting the region into the box $(0,\delta) \times (r - \delta, r + \delta)$. In the new coordinates,
the finite Taylor expansion $R(x,y) = \sum_{a,b < K}R_{ab}x^ay^b + O(|x|^K + |y|^K)$ becomes of the form
$$R(x,x^My') = x^{\alpha }R_M(1,y') + x^{\alpha + \zeta}s(x,y') + O(|x|^K + |x|^{KM}|y'|^K) \eqno (2.11)$$
Here $\zeta > 0$ and $f(x,y')$ is a polynomial in $y'$ and a fractional power of $x$. For any $\epsilon' > 0$,
if $\delta$ is sufficiently small we have $|R_M(1,y') - R_M(1,r)| < \epsilon'$ for all $|y' - r| 
< \delta$. Equivalently, $|R_M(1,y') - R_M(1,r)|x^{\alpha} < \epsilon'x^{\alpha}$.
Furthermore, if $\delta$ is sufficiently small $x^{\alpha + \zeta}|s(x,y')|$ and the 
$O(|x|^K + |x|^{KM}|y'|^K)$ term are less than $\epsilon' x^{\alpha}$ whenever $x$ and $|y' - r|$
are sufficiently smalll. Combining, if $\delta$ is sufficiently small then on our domain we have
$$|R(x,x^My') -  x^{\alpha}R_M(1,r)| < 3 \epsilon' x^{\alpha}R_M(1,r)\eqno (2.12)$$
Translating this back into the original coordinates, we have
$$|R(x,y) - R_M(x,y)| < 3\epsilon' x^{\alpha} R_M(1,r) \eqno (2.13)$$
Taking $\epsilon = 3\epsilon' R_M(1,r)$ gives us the lemma and we are done.

We now are in a position to set up the proof of the main theorem, Theorem 1.1. Recall we have a surface
$Q$ with a distiguished point $(x_0,y_0,z_0)$ that is the graph of some smooth function $g(x,y)$ defined 
near $(x_0,y_0)$ such that $\nabla g(0,0) = (0,0)$. Suppose the assumptions of Theorem 1.1 hold. Then, 
after a linear coordinate change if necessary, we may assume $f(x,y) = g(x_0 + x,y_0 + y) - z_0$ has a 
generic adapted coordinate system on a neighborhood of the origin for which the exceptional situation b) of 
Theorem 1.1 does not occur. Therefore there is a smooth $\psi(x)$
with $\psi(0) = 0$ such that $F(x,y) = f(x,y + \psi(x))$ is in generic adapted coordinates, and if 
$N(F)$ intersects the bisectrix in the interior of a compact edge $e$ then $e$ does not satisfy the 
exceptional situation b) of the statement of Theorem 1.1. 

\noindent {\bf Definition of slivers for $F(x,y)$}

We now use Lemmas 2.2 and 2.3 on $F(x,y)$ and its various $y$ derivatives to subdivide a small neighborhood $B$
of $(0,0)$ into "slivers" containing the origin. The case where $N(F)$ has exactly one vertex (which is 
therefore on the $y$-axis) is easier and will be treated separately, so in the following we always assume
$N(F)$ contains multiple vertices. Denote the vertices of $N(F)$ above the
bisectrix by $v_1,...,v_k$ where if $i < j$ then $v_i$ is below $v_j$. Let $e_i$ denote the edge of $N(F)$ whose 
upper vertex is $v_i$; if $v_1$ is the lowest vertex of $N(F)$ then we just do not define $e_1$.
Let $m_i = -{1 \over M_i}$ denote the slope of $e_i$. Write $v_i = (a_i,b_i)$; observe $b_i \geq 2$ for 
all $i$ since $(a_i,b_i)$ lies above the bisectrix and we are assuming $F(0,0) = 0$ and $\nabla F(0,0) = 
(0,0)$. 

For $v_i$ that is the intersection of two compact edges of $N(F)$, define $D_i$ to be the set $\{(x,y)
\in B:  N_0 |x|^{M_i} < |y| < {1 \over N_0} |x|^{M_{i+1}}\}$. Here $N_0$ is large enough so that we may
invoke Lemma 2.2 and say there are  $c_0,c_1 > 0$ such that (assuming $B$ is small enough) for 
$m = 0,1,2$ on $D_i$ we have
$$ c_1 |x|^{a_i} |y|^{b_i - m} > |{\partial^m F \over \partial y^m}(x,y)| > c_0 |x|^{a_i} |y|^{b_i - m}
\eqno (2.14)$$
It should be pointed out that if $b_i = m = 2$, then $(2.14)$ holds by applying Lemma 2.3, reversing the roles
of the $x$ and $y$ axes and setting $r = 0$. If $v_i$ is the upper vertex of $N(F)$, we define
$D_i$ to be the points where $N_0 |x|^M_i < |y|$, in which case $(2.14)$ 
still holds on $D_i$ by Lemma 2.2.

We next subdivide the set $B - \cup_i D_i$ into some 
slivers touching the origin amenable to the analysis of this paper. We only describe the slivers for $x > 0$;
the ones where $x < 0$ are defined analogously. Note that the points of $B - \cup_i D_i$ where $x > 0$ can 
be written as $\cup_i C_i$, where 
$$C_1 = \{(x,y) \in B: x > 0,\,\,|y|  < N_0 x^{M_1} \} \eqno (2.15a)$$
$$C_i = \{(x,y) \in B:  x > 0,\,\,  {1 \over N_0} x^{M_i} < |y|  < N_0 x^{M_i}\} \,\,\,\,\,\,\,\,\,(i > 1) \eqno (2.15b)$$
(In the event that $v_1$ is the lowest vertex of $N(F)$, one takes $M_2$ in
$(2.15a)$ and then $(2.15b)$ is valid for $i > 2$).
Suppose $r$ is such that $F_{e_i}(1,r) \neq 0$, but $\partial_y F_{e_i}(1,r)$ has a zero of order 
greater than $d^* - 1  = \max(2, d(F)) - 1$ at $r$. Then if $e_i$ intersects the bisectrix in its interior,
by assumption the exceptional case of Theorem 1.1 part b does not occur. If $e_i$ does not intersect the 
bisectrix in its interior, then by Lemma 2.1 $e_i$ has slope -1 and intersects the $y$-axis. In either event, if $i = 1$ and $|r| < N_0$ or
$i > 1$ and ${1 \over N_0} < |r| < N_0$, we define $E_{ir}$ to be the sliver 
$$E_{ir} = \{(x,y) \in B:  x > 0, \,\, (r - \delta_r) x^{M_i} < y  < (r + \delta_r) x^{M_i} \}\eqno (2.16)$$
Here $\delta_r$ is a small constant to be determined by our future arguments. We will refer to the 
(finitely many) $E_{ir}$ occurring as $E_{ij}$ in the rest of this paper.

For any $r$ other than these,
we may let $1 \leq k \leq d^*$ and $\delta_r, C_r > 0$ be such that 
on $[r - \delta_r, r + \delta_r]$ we have $|\partial_y^{k} F_{e_i}(1,y)| > C_r$. If $F_{e_i}(1,r) \neq 0$ 
this follows from the above definition of the $D_i$ and if $F_{e_i}(1,r) = 0$ it follows from the fact
that any zero of any $F_{e_i}(1,y)$ has order at most $d(F)$ in adapted coordinates. As a result, on the set 
$B_r = \{(x,y) \in B: x > 0,\,\, (r - \delta_r) x^{M_i} \leq y  \leq (r + \delta_r) x^{M_i}\}$, given that
$e_i$ contains $(a_i,b_i)$ we have
$$|\partial_y^{k} F_{e_i}(x,y)| > C_r x^{a_i+ M_i(b_i - {k})} \eqno (2.17)$$
By applying Lemma 2.3 to $\partial_y^k F(x,y)$ we can assume $B$ is small enough that we also have
$$|\partial_y^{k} F(x,y) - \partial_y^k F_{e_i}(x,y)| < {C_r \over 2}x^{a_i+ M_i(b_i - {k})} \eqno (2.18)$$
Putting $(2.17)$ and $(2.18)$ together, on $B_r$ we have
$$ |\partial_y^{k} F(x,y)| > {C_r \over 2} x^{a_i+ M_i(b_i - {k})} \eqno (2.19)$$
By compactness, we can write $B - \cup_i D_i - \cup_{ij} E_{ij}$ as the union of finitely many slivers on which 
$(2.19)$ is satisfied. For a given edge $e_i$, we write the slivers for which $k = 1$
as $F_{ij}$ and the slivers for which $k > 1$ by $G_{ij}$. We denote the value
of $k$ corresponding to a given $G_{ij}$ by $k_{ij}$. Note that each $k_{ij} \leq d^*$.

The above decompositions were for the case where $N(F)$ had more than one vertex. When $N(F)$ has exactly
one vertex, since $F$ it is in generic adapted coordinates it is of the form $(0,k)$. In this case, we 
simply designate a neighborhood of the origin as a single $G_{ij}$, with $k_{ij} = k$. In general, the
arguments for this $G_{ij}$ will be simplified versions of the $G_{ij}$ arguments for the multivertex case.

\noindent Let $v(F)$ denote the set of vertices of $N(F)$, and define $F^*(x,y)$ by 
$$F^*(x,y)= \big(\sum_{(v_1,v_2) \in v(F)} (x^{v_1}y^{v_2})^2\big)^{1 \over 2}$$
The function $F^*(x,y)$ will be used in defining the damping function.
To this end, first note that $N((F^*)^2)$ is the double of $N(F)$ and therefore $d((F^*)^2) = 2d(F)$.
As a result, by [V], 
$F^*(x,y)^{-t}$ is integrable on a neighborhood of the origin iff $t < {1 \over 2}d((F^*)^2) = d(F)$. We apply Lemmas 2.2 and 2.3 
to $F^*(x,y)^2$ in place of $F(x,y)$ and obtain that if the $N_0$ in the definition of $D_i$ were chosen
sufficiently large, then there is a constant $C_0 > 0$ such that on each $D_i$ we have
$${1 \over C_0}|x^{a_i}y^{b_i}| <  F^*(x,y) < C_0|x^{a_i}y^{b_i}| \eqno (2.20)$$
Similarly, by Lemma 2.3 on each $E_{ij}$, $F_{ij}$ and $G_{ij}$, the constant $C_0$ can be taken so that we have
$${1 \over C_0}|x|^{a_i + M_ib_i} < F^*(x,y) < C_0 |x|^{a_i + M_ib_i} \eqno (2.21)$$
We now subdivide the surface $Q$ near our distinguished point $q_0$ in accordance with the above
subdivisions, applied to the function $f(x,y)$ in generic adapted coordinates. In other words, we let
$\psi(x)$ be such that $F(x,y) = f(x, y + \psi(x))$ is in generic adapted coordinates such that the 
exceptional cases of Theorem 1.1 do not hold, and define $D_i$, $E_{ij}$,
$F_{ij}$, and $G_{ij}$ to be the above slivers as defined for $F(x,y)$. We next transfer these slivers into
the original coordinates of the surface $Q$; let $D_i'$ be the portion of $Q$ above the set
$\{(x,y): (x - x_0,y - y_0 - \psi(x - x_0)) \in D_i\}$, with the analogous definitions for $E_{ij}'$,
$F_{ij}'$, and $G_{ij}'$. We also will have use for $D_i$ in the original nonadapted 
coordinates, centered at $(x_0,y_0)$. To that end we let $D_i'' = \{(x,y): (x, y - \psi(x)) \in D_i\}$,
making the analogous definitions for $E_{ij}''$, $F_{ij}''$, and $G_{ij}''$.  

\noindent The next lemma will be
useful in bounding the contribution of our integrals over $E_{ij}''$ in the $L^2$ estimates of section 4. 

\noindent {\bf Lemma 2.4}. There is a constant $C$ such that on $E_{ij}''$ we have
$|\partial_{xx} f(x, y)| \geq C|x|^{a_i + b_iM_i - 2}$.

\noindent {\bf Proof.} We consider slivers for which $x > 0$ as the $x < 0$ slivers are entirely analogous.
Recall each $E_{ij}$ is a region of the form $\{(x,y): 0 < x < \mu, (r - \nu)x^{M_i} < y < (r + \nu)x^{M_i}\}$,
where $F_{e_i}(1,r) \neq 0$, but where $\partial_y^lF_{e_i}(1,r) = 0$ for at least $l = 1, 2$.
Under the map $(x,y) \rightarrow (x, y - rx^{M_i})$, the set 
$E_{ij}$ becomes the region $E_{ij}'''$ = $\{(x,y): 0 < x < \mu, |y| < \nu x^{M_i}\}$, and if $G(x, y)$ 
denotes $F(x,y + rx^{M_i})$ then  $G_{e_i}(1,0) \neq 0$ but
$\partial_y^lG_{e_i}(1,0) = 0$ for $l = 1, 2$. 

In terms of Newton polygons, the above can be translated as follows. Since $e_i$ is an edge of 
$N(F)$ with equation $x + M_i y = a_i + M_ib_i$, $N(G)$ has an edge with the same equation which goes all
the way to the $x$-axis since $G_{e_i}(1,0) \neq 0$.
As a result, for $l = 1,2$ $N(\partial_x^l G)$ has an edge with equation $x + M_i y = a_i + M_ib_i - l$ 
which extends to the $x$ axis. On the other hand, since $\partial_y G_{e_i}(1,y)$ has a zero of order at least
two at $0$, for $l = 1,2$ $N(\partial_y^l G)$ intersects the line with equation $x + M_i y = a_i + 
(M_i - l)b_i$ but does not contain $(a_i + (M_i - l)b_i,0)$. Using this fact for $l = 1$ and taking an $x$ 
derivative shows that $N(\partial_{xy} G)$ intersects the line with equation $x + M_i y = a_i - 1 + 
(M_i - 1)b_i$ but does not contain $(a_i - 1 + (M_i - 1)b_i,0)$.
Using $(2.11)$ in conjunction with these latter observations concerning the 
Newton polygons, we obtain that for any $\eta > 0$, if $\nu$ were chosen sufficiently
small, then on $E_{ij}'''$ for $l = 1, 2$ we have
$$|\partial_y^l G(x,y)| < \eta |x|^{a_i + (M_i - l)b_i} \eqno (2.22a)$$
$$|\partial_{xy} G(x,y)| < \eta |x|^{a_i - 1 + (M_i - 1)b_i} \eqno (2.22b)$$
On the other hand Lemma 2.3 in conjuction with the above conditions on the Newton polygons of the 
$x$-derivatives ensures that for some $c_1 > 0$, on $E_{ij}'''$ for $l = 1, 2$ we have
$$|\partial_x^l G(x,y)| \geq c_1 |x|^{a_i + M_ib_i - l} \eqno (2.22c)$$
Next, as mentioned after $(2.15)$, $E_{ij}$ is only defined in two situations.
The first is when $e_i$ intersects the bisectrix in its interior and the function $\psi(x)$ such that $F(x,y) = 
f(x, y + \psi(x))$ has a zero of order at least $M_i$ at $x = 0$. In the second situation, $e_i$ does not intersect the
bisectrix in its interior, but by Lemma 2.1 $M_i = 1$ and thus $\psi(x)$ still has a zero of order at least
$M_i$ at $x = 0$. In either case, since $f(x,y) = F(x, y - \psi(x)) = G(x, y - \psi(x) + rx^{M_i})$, we can
write $f(x,y) = G(x,y + \xi(x))$ where $\xi(x)$ has a zero of order at least $M_i$ at $x = 0$. Applying the
chain rule, we get that on $E_{ij}''$ we have
$$\partial_{xx} f(x,y) = \partial_{xx}G(x, y + \xi(x)) + 2\xi'(x)\partial_{xy} G(x, y + \xi(x))$$
$$+ \xi''(x)\partial_yG(x, y + \xi(x)) + (\xi'(x))^2\partial_{yy}G(x, y + \xi(x))$$
Equation $(2.22c)$ ensures that $|\partial_{xx}G(x, y + \xi(x))| \geq c_1 |x|^{a_i + M_ib_i - 2}$, and
equations $(2.22a)-(2.22b)$ coupled with the fact that $\xi(x)$ has a zero of order at least $M_i$ at $x = 0$ 
ensure that the remaining terms can be made less than any $\eta'
|x|^{a_i + M_ib_i - 2}$. We conclude that $|\partial_{xx} f(x,y)|$ is at least ${c_1 \over 2}
|x|^{a_i + M_ib_i - 2}$ on $E_{ij}''$ if $\nu$ were chosen appropriately small. This completes the proof of 
the lemma.

\noindent {\bf Definition of the damping factor.}

We now define the damping factor on the surface $Q$ in a neighborhood of $q_0 = (x_0,y_0, g(x_0,y_0))$.
Above a point $(x,y)$ it will be of the form $e^{z^2}|h(x,y)|^z|D(x,y)|^{\delta z}$, where $D(x,y)$ is 
the Hessian determiant of $g$ at $(x,y)$ and $\delta$ is a small positive number to be determined by our 
arguments. The function $h(x,y)$ will be defined in the form $H(x - x_0 ,y - y_0 - \psi(x - x_0))$, where
$H(x,y)$ is expressed in terms of $F(x,y) = f(x,y + \psi(x))$. To this end, once again
let $d^* = \max(2,d(F))$. Since $F(x,y)$ is in adapted coordinates, $d(F) = h(q_0)$ and thus equivalently
we have $d^* = \max(2, h(q_0))$. 

On each $D_i$, $E_{ij}$, $F_{ij}$, as well as on each $G_{ij}$ with $k_{ij} = 2$, we define $H(x,y)$
to be $F^*(x,y)^{{1 \over 2} - {1 \over d^*}}$. On the remaining $G_{ij}$, if $N(F)$ has just one vertex
we let $H(x,y) = |{\partial^2F \over \partial y^2}(x,y)|^{1 \over 2}$, while if $N(F)$ has
multiple vertices let $H(x,y) = 
|x|^{M_i - {a_i + M_ib_i \over d^*}}|{\partial^2F \over \partial y^2}(x,y)|^{1 \over 2}$. One thing worth 
mentioning concerning these latter $H(x,y)$ is the following. Since $|y| < C|x|^{M_i}$ on $G_{ij}$, by
$(2.11)$ with $\alpha = a_i + M_ib_i$, on $G_{ij}$ one has
$|{\partial^2F \over \partial y^2}(x,y)| < C|x|^{a_i+ M_i(b_i - 2)}$. As a result,
$|H(x,y)|\leq C |x|^{(a_i+ M_ib_i)({1 \over 2} - {1 \over d^*})}$ and in view of $(2.21)$ this gives
$$|H(x,y)|\leq C F^*(x,y)^{{1 \over 2} - {1 \over d^*}}\eqno (2.23)$$
Note the right-hand side of $(2.23)$ is exactly $CH(x,y)$ for the other regions.

\noindent {\bf 3. $L^{\infty}$ estimates}

Define the operator $M_z$ to be the maximal operator $(1.1)$ with respect to the measure 
$e^{z^2}|h(x,y)|^z |D(x,y)|^{\delta z}\phi(q)d\sigma(q)$ in place of $\phi(q)d\sigma(q)$. Note that for 
$z = 0$, $M_z$ is exactly $M$. These will be the analytic family of maximal operators used in proving Theorem
1.1 as described at the end of section 1. We first prove the $L^{\infty}$ to $L^{\infty}$ boundedness 
properties of the $M_z$ we need.

\noindent {\bf Theorem 3.1.} Write $d = d(F)$. If $d > 2$, then for any $s > -{2 \over d - 2}$, if $\delta$ is sufficiently
small (depending on $s$) then there exists a constant $C$ such that $||M_z||_{L^{\infty} \rightarrow 
L^{\infty}} < C$ for all $z$ with $Re(z) = s$. If $d \leq 2$, the same holds for any $s \in \R$.

\noindent {\bf Proof.} We will use the fact that $||M_z||_{L^{\infty} \rightarrow L^{\infty}}$ is
bounded by the $L^1$ norm of the damping function of $M_z$. We consider the case $d = d(F) \leq 2$ first. 
Note that for each $D_i$, $E_{ij}$ or $F_{ij}$ the function $h(x,y)$ is
given by $F^*(x - x_0, y - y_0 - \psi(x-x_0))^{{1 \over 2} - {1 \over d^*}} = 1$. Since in adapted 
coordinates each $k_{ij}$ is at most $d^* \leq 2$ here, for each $G_{ij}$ that may appear $h(x,y)$ is also 
just 1. We conclude that the
damping factor always equal to $e^{z^2}|D(x,y)|^{-\delta z}$. On a given vertical line $Re(z) = s$, this
has magnitude bounded by $C_s|D(x,y)|^{-\delta s}$.  Since $D(x,y)$ is 
assumed to be of finite type at the origin, there is some $\epsilon > 0$ for which
$|D(x,y)|^{-\epsilon}$ is integrable on a neighborhood of the origin. Hence as long as $\delta < {\epsilon
\over s}$, the damping factor is integrable, with integral uniformly bounded on $Re(z) = s$. This is 
exactly what we needed to prove.

Now suppose $d > 2$. On any $D_i'$, $E_{ij}'$, $F_{ij}'$, or $G_{ij}'$ with $k_{ij} = 2$, 
$h(x,y) = F^*(x - x_0, y - y_0 - 
\psi(x-x_0))$ and thus the damping factor has magnitude $F^*(x - x_0, y - y_0 - \psi(x-x_0))^{({1 \over 2} - 
{1 \over d})Re(z)}|D(x,y)|^{\delta Re(z)}$. As mentioned at the end of section 2, a result of Varchenko
says that $F^*(x,y)^{-t}$ is integrable on a neighborhood of the origin iff $t < {1 \over d}$. Thus
the same is true for $F^*(x - x_0, y - y_0 - \psi(x-x_0))^{-t}$.
As a result, $F^*(x - x_0, y - y_0 - \psi(x-x_0))^{({1 \over 2} - {1 \over d})Re(z)}$ is integrable on a 
neighborhood of the origin iff $Re(z) > {-{2 \over d - 2}}$. Consequently, by Holder's inequality,
for fixed $s > {-{2 \over d - 2}}$, by choosing $\delta$ sufficiently small we have that when $Re(z) = s$ 
the damping factor is integrable over any $D_i'$, $E_{ij}'$, $F_{ij}'$, or $G_{ij}'$ with $k_{ij} =2$, 
with integral uniformly bounded in $Im(z)$. This is what we need here.

We now move on to the $G_{ij}'$ with $k_{ij} > 2$. If $N(F)$ has one vertex, we use the fact that 
${\partial^2F \over \partial y^2}(x,y)$ has nonvanishing $(k_{ij} - 2)$th derivative in the $y$ direction,
where $k_{ij} \leq d$. Thus if $Re(z) =  s >  
-{2 \over d - 2}$, then  $|H(x,y)|^z = |{\partial^2 F\over \partial y^2}(x,y)|^{z \over 2}$ is integrable
in $y$ with integral uniformly bounded in $Im(z)$. Making $\delta$ sufficiently small and using Holder's
inequality again gives the desired result. Suppose now $N(F)$ has multiple 
vertices. We will consider those $G_{ij}$ for which $x > 0$ as the $G_{ij}$ for which $x < 0$ are done in the 
same way. In the $G_{ij}$ coordinates we can write the damping function as
$$|\tilde{D}(x,y)|^{\delta z} x^{M_i - {a_i + M_ib_i \over d}}|{\partial^2F \over \partial y^2}(x,y)|^{1 \over
2}$$
Recall that $G_{ij}$ is of the form $\{(x,y): 0 < x < \eta, (r - \delta)x^{M_i} < y < (r + \delta)x^{M_i}\}$.
Analogous to $(2.11)$, on the box
$(0, \eta) \times (r - \delta, r + \delta)$ for any $K$ one can write 
$$F(x,x^{M_i}Y) = x^{\alpha_i}F_{M_i}(1,Y) + x^{\alpha_i + \zeta}s(x,Y) + O(|x|^K + |x|^{KM_i}|Y|^K) 
\eqno (3.1a)$$
Here $s(x,Y)$ is a polynomial in $Y$ and a fractional power of $x$. Analogous expressions hold for various
$y$ derivatives of $F$; for example, we can write
$${\partial^2F \over \partial y^2} (x,x^{M_i}Y) = x^{\alpha_i - 2M_i}\partial_{yy} F_{M_i}(1,Y) + 
x^{\alpha_i - 2M_i + \zeta}\tilde{s}(x,Y) + O(|x|^K + |x|^{KM_i}|Y|^K) \eqno (3.1b)$$
The constant $\delta$ was chosen small enough that $|\partial_Y^{k_{ij}}F_{M_i}(1,Y)| > C$ on $(r - \delta, 
r + \delta)$ for some positive $C$, where $k_{ij} \leq d$. As a result, shrinking $\eta$ if necessary we
can assume that on $(0, \eta) \times (r - \delta, r + \delta)$ we have ${\partial^2F \over \partial y^2} 
(x,x^{M_i}Y)= x^{\alpha_i - 2M_i}S(x,Y)$ where 
$$|\partial_Y^{k_{ij} - 2}S(x,Y)| > C \eqno (3.2)$$
We now let $x = X^{1 \over M_i + 1}$, so that ${\partial^2F \over \partial y^2} (X^{1 \over M_i + 1},X^{M_i
\over M_i + 1}Y)$ is of the form $X^{\alpha_i - 2M_i \over M_i + 1}T(X,Y)$ where
$$|\partial_Y^{k_{ij} - 2}T(X,Y)| > C \eqno (3.2')$$
We make this coordinate change so that the change from $(x,y)$ coordinates to $(X,Y)$ has constant Jacobian
determinant. This ensures that a power of the damping function is integrable in the $(x,y)$ coordinates iff
it is integrable in the $(X,Y)$ coordinates. In the old coordinates, the damping function is 
$|\tilde{D}(x,y)|^{\delta z}$ times  $x^{({dM_i - a_i - M_ib_i \over d})z}|{\partial^2F \over \partial y^2}(x,y)|^
{z \over 2}$ so in the new coordinates it is of the form 
$$|\bar{D}(X,Y)|^{\delta z} \times X^{{dM_i - a_i - M_ib_i \over d(1 + M_i)}z}
|X^{\alpha_i - 2M_i \over M_i + 1}{\partial^2T \over \partial Y^2}(X,Y)|^{z \over 2} $$
$$= |\bar{D}(X,Y)|^{\delta z} \times X^{{d{\alpha_i \over 2} - a_i - M_ib_i \over d(M_i + 1)}z}
|{\partial^2T \over \partial Y^2}(X,Y)|^{z \over 2} \eqno (3.3)$$
Note that by $(3.2')$ the function ${\partial^2T \over \partial Y^2}(X,Y)$ has $(k_{ij} - 2)$th 
derivative uniformly bounded below in the $y$ variable, for fixed $x$, where $k_{ij} \leq d$. As a result,
if $Re(z) =  s >  -{2 \over d - 2}$, then  $|{\partial^2 T\over \partial Y^2}(X,Y)|^{z \over 2}$ is integrable in $y$
with integral uniformly bounded in $Im(z)$.  Thus for such $z$ we have
$$\int_{r - \delta}^{r + \delta}\int_0^{\eta} |x^{\big({ d{\alpha_i\over 2} - a_i - M_ib_i \over d(M_i + 1)}
\big)}|^{z \over 2} |{\partial^2T \over \partial Y^2}(X,Y)|^{z \over 2} \,dx\,dy< 
C\int_0^{\eta} x^{s{d{\alpha_i \over 2}- a_i - M_ib_i \over d(M_i + 1)}}\,dx\eqno (3.4)$$
This will be finite if the exponent of $x$ is greater than $-1$. Substituting $\alpha_i = a_i + M_ib_i$, the 
exponent in $(3.4)$ is $s{d - 2 \over 2d} {a_i + M_i b_i \over 1 + M_i}$. 
Note that ${a_i + M_ib_i \over 1 + M_i}$ is the $x$-coordinate of the intersection of the bisectrix with 
the edge $e_i$, which is at most $d$. Hence ${d - 2 \over 2d} {a_i + M_i b_i \over 1 + M_i} < 
{d - 2 \over 2}$. Therefore if $s > {2 \over d - 2}$, the exponent in $(3.4)$ is greater
than $-1$ and thus the right-hand factor of $(3.3)$ 
is integrable for $Re(z) = s$, uniformly in $Im(z)$. As in the previous argument, by
making $\delta$ in the $|\bar{D}(X,Y)|^{\delta z}$ factor sufficiently small by Holder's inequality the same 
will be true for the entire damping factor $(3.3)$. This what we needed to prove and we are done.

\noindent {\bf 4. $L^2$ estimates.}

We now move to proving the $L^2$ bounds needed in the proof of Theorem 1.1. As indicated in section 1, we
will be utilizing Theorem 1.2 that follows from [SoSt]. Letting $\phi(q)d\sigma(q)$ be the surface measure
of $(1.1)$, we define the measure $\sigma_z$ by 
$d\sigma_z(q) = e^{z^2}|h(x,y)|^z|D(x,y)|^{\delta z}\phi(q)d\sigma(q)$. 
Since we will using Theorem 1.2, we examine its Fourier transform $\hat{\sigma_z}(\lambda)$, given by
$$\hat{\sigma_z}(\lambda) = e^{z^2} \int e^{-i\lambda_1 g(x,y) - i\lambda_2x - i\lambda_3 y}|h(x,y)|^z 
|D(x,y)|^{\delta z}\phi(x,y) \,dx\,dy \eqno (4.1)$$
Here $\phi(x,y)$ denotes some cutoff function on a neighborhood of the origin. We always shift by 
$(x_0,y_0)$, so that our integrals are over a small neighborhood of the origin. Thus up
to an ignorable factor of magnitude 1, $(4.1)$ is given by
$$\hat{\sigma_z}(\lambda) = e^{z^2} \int e^{-i\lambda_1 f(x,y) - i\lambda_2x - i\lambda_3 y}|H(x,y - \psi(x))|^z 
|D(x,y)|^{\delta z}\phi^*(x,y) \,dx\,dy \eqno (4.2)$$
In $(4.2)$, $\psi(x)$ is the function taking $f(x,y)$ into its adapted coordinates and $D(x,y)$ now denotes
the Hessian determinant of $f$ at $(x,y)$.  For the analysis of the part of $(4.2)$ coming from the $D_i''$,
$F_{ij}''$, and $G_{ij}''$, we 
will transfer into the adapted coordinates of $f(x,y)$ and use Van der Corput-type arguments in the 
$y$-variable. For the $E_{ij}''$ we will remain in the original coordinates, and use a Van 
der Corput-type argument in the $x$-variable in conjunction with Lemma 2.4.

We will prove that the conditions of Theorem 1.2 hold by virtue of the following theorem, whose proof will
comprise most of the rest of this paper.

\noindent {\bf Theorem 4.1.} Suppose $s > 1$. Then if the constant $\delta$ used in the exponent of 
$|D(x,y)|$ is sufficiently small, then there are constants $C, \epsilon$ independent of $Im(z)$ such that 
if $Re(z) = s$ then for any multiindex $\alpha$ with $|\alpha| = 0,1$ we have
$$|\partial^{\alpha} \hat{\sigma_z} (\lambda)| < C(1 + |\lambda|)^{-{1 \over 2} - \epsilon} \eqno (4.3)$$

\noindent {\bf Proof.} We will only prove $(4.3)$ for $|\alpha| = 0$ as the $|\alpha| = 1$ cases are 
identical other than having a different cutoff function $\phi(x,y)$. Recall that we are assuming that $f(0,0) =
0$ and $\nabla f(0,0) = 0$. So if  $|\lambda_2|$ or $|\lambda_3|$ is the maximal $|\lambda_i|$
one may integrate by parts in $x$ or $y$ respectively and get that $|\hat{\sigma}(\lambda)| < C|\lambda|^{-1}$,
which is better than the estimate that we need. Hence for the remainder of this paper we will always assume 
that $|\lambda_1|$ is at least as large as $|\lambda_2|$ and $|\lambda_3|$.

Let $\alpha(x)$ be an even function on $\R$ that is equal to 1 for $|x| \leq {1 \over 2}$, zero for $|x| > 1$, and 
is monotone decreasing on $\R^+$. Let $\beta(x) = 1 - \alpha(x)$. For constants $\delta_1$ and $N_1$ to be 
determined by our arguments, we express $(4.2)$ as $I_1 + I_2 + I_3$, where
$$I_1(\lambda) =  e^{z^2} \int e^{-i\lambda_1 f(x,y) - i\lambda_2x - i\lambda_3 y}|H(x,y - \psi(x))|^z
|D(x,y)|^{\delta z} $$
$$\times \alpha(|\lambda|^{N_1} D(x,y)) \phi^*(x,y) \,dx\,dy \eqno (4.4a)$$
$$I_2(\lambda) =  e^{z^2} \int e^{-i\lambda_1 f(x,y) - i\lambda_2x - i\lambda_3 y}|H(x,y- \psi(x))|^z
|D(x,y)|^{\delta z}$$
$$\times \big( \alpha(|\lambda|^{\delta_1} D(x,y) - \alpha(|\lambda|^{N_1} D(x,y))\big)
\phi^*(x,y) \,dx\,dy \eqno (4.4b)$$
$$I_3(\lambda) =  e^{z^2} \int e^{-i\lambda_1 f(x,y) - i\lambda_2x - i\lambda_3 y}|H(x,y- \psi(x))|^z
|D(x,y)|^{\delta z}$$
$$\times \beta(|\lambda|^{\delta_1} D(x,y)) \phi^*(x,y) \,dx\,dy \eqno (4.4c)$$

The analysis of $I_2(\lambda)$ will be the crux of the argument. The contribution to $(4.3)$ due to 
$I_1(\lambda)$ is easily shown to decrease rapidly in $|\lambda|$. Specifically, since $D(x,y)$ is being
assumed to be of finite-type in
a neighborhood of the origin, if $N_1$ is large enough the measure of the points where $|D(x,y)| < 
|\lambda|^{-{N_1}}$ will be less than ${1 \over |\lambda|}$. As a result, the integrand of $I_1(\lambda)$
is nonzero on a set of measure at most ${1 \over |\lambda|}$. Since all the factors in $(4.4a)$ are 
uniformly bounded on a line $Re(z) = s$ with $s > 1$, this gives that $|I_1(\lambda)| < C|\lambda|^{-1}$, better than
what is needed.

\noindent {\bf Bounding $|I_3(\lambda)|$}.

Note that on the support of the integrand
of $I_3(\lambda)$, the Hessian determinant $D(x,y)$ is at least ${1 \over 2}|\lambda|^{-\delta_1}$. The
idea is that if $\delta_1$ were actually zero, then on this support the Hessian would be bounded below and
we would get an estimate $|I_3(\lambda)| < C|\lambda|^{-1}$. Although $\delta_1$ is not zero, if it is 
sufficiently small we still get an estimate $|I_3(\lambda)| < C|\lambda|^{-{1 - t}}$ for any given small 
$t$, an estimate better than what is needed.

We proceed as follows. For a sufficiently small $c > 0$ to be determined by our arguments, we divide the 
support of the integrand of $I_3(\lambda)$ into squares of diameter $c|\lambda|^{-\delta_1}$. We will show
that the contribution to $I_3(\lambda)$ from each such square is at most $c|\lambda|^{-{3 \over 5}}$
if $\delta_1$ is sufficiently small. Adding this over all these squares, this gives an estimate better 
than needed. 

Let $S$ be any such square. Since $D(x,y)$ is of finite type, we may let $u$ and $v$ be nonparallel directions
such that for some $k$, $\partial_u^kD(x,y)$, $\partial_v^kD(x,y)$, and $\partial_u\partial_v^{k-1}D(x,y)$ 
are nonvanishing on the support
of the integrand of $I_3(\lambda)$. We can similarly assume that there are $k',k''$ such that 
$\partial_u^{k'}\big(F^*(x,y - \psi(x))^2\big)$, $\partial_v^{k'}\big(F^*(x,y - \psi(x))^2\big)$, $\partial_u^
{k''}f(x,y)$, and $\partial_u^{k''}f(x,y)$ are nonvanishing on any $S$.  Let $a_1 = 
-\partial_u ({\lambda_2\over \lambda_1}x + {\lambda_3\over \lambda_1}y)$ and $a_2 = -\partial_v 
({\lambda_2\over \lambda_1}x + {\lambda_3\over \lambda_1}y)$. Note $a_1$ and $a_2$ are constants. Define the
sets $S_1$, $S_2$, and $S_3$ by
$$S_1 = \{(x,y) \in S: |\partial_uf(x,y) - a_1| > |\lambda|^{-{1 \over 3}}\} \eqno (4.5a)$$
$$S_2 = \{(x,y) \in S: |\partial_uf(x,y) - a_1| \leq |\lambda|^{-{1 \over 3}}, |\partial_vf(x,y) - a_2| >
|\lambda|^{-{1 \over 3}} \} \eqno (4.5b)$$
$$S_3 = \{(x,y) \in S: |\partial_uf(x,y) - a_1| \leq |\lambda|^{-{1 \over 3}}, |\partial_vf(x,y) - a_2| \leq
|\lambda|^{-{1 \over 3}} \} \eqno (4.5c)$$
Correspondingly, write the contributions to $I_3(\lambda)$ from $S_1$, $S_2$, and $S_3$ as $J_1(\lambda)$, $J_2(\lambda)$, and
$J_3(\lambda)$ respectively. To analyze $J_1(\lambda)$, integrate the integrand of $(4.4c)$ by parts in the $u$ direction,
integrating $\lambda_1(\partial_uf(x,y) - a_1)e^{-i\lambda_1 f(x,y) - i\lambda_2x - i\lambda_3 y}$ in and differentiating
${1 \over \lambda_1(\partial_uf(x,y) - a_1)}$ times the remainder of the integrand. We get several terms depending
on where the derivative lands. If it lands on the $\phi(x,y)$ factor, then each factor in the term is 
bounded above by a constant, with the exception of the ${1 \over \lambda_1(\partial_uf(x,y) - a_1)}$ factor,
which is bounded in absolute value by $|\lambda_1|^{-{2 \over 3}}$ on $S_1$. Hence $|J_1(\lambda)|$ is
at most $C|\lambda_1|^{-{2 \over 3} + 2\delta_1}$, which is bounded by $C|\lambda|^{-{3 \over 5}}$, the desired estimate.

Next, we consider the term where the derivative lands on the ${1 \over \lambda_1(\partial_uf(x,y) - a_1)}$
factor. We take absolute values of the entire integrand and bound it above by ${C \over |\lambda|}
{\partial_u^2 f(x,y) \over (\partial_uf(x,y) - a_1)^2}$. We integrate this in the $u$ direction as in the
proof of the Van der Corput lemma; the assumed condition that $|\partial_u^{k''}f(x,y)|$ is bounded below
ensures that we integrate over boundedly many intervals on which $\partial_uf(x,y) - a_1$ is monotone and
thus ${\partial_u^2 f(x,y) \over (\partial_uf(x,y) - a_1)^2}$ integrates back to ${1 \over (\partial_uf
(x,y) - a_1)}$. We end out with a bound of $C|\lambda_1|^{-{2 \over 3}  + \delta_1} < C|\lambda_1|^
{-{3 \over 5}}$. 

If the derivative lands on the $|D(x,y)|^{\delta z}$ factor, we argue similarly.
We take absolute values and integrate in the $u$ direction, this time using that $|\partial_u^kD(x,y)|$ is
bounded below on the integrand to ensure that there are boundedly many intervals on which $D(x,y)$ is 
monotone and thus on which we can integrate back its $u$-derivative. In the (extremely
rare) case that only $k = 0$ can be used, $|D(x,y)|^{\delta z}$ is a smooth function and the term behaves
as in the case where the derivative lands on $\phi(x,y)$. If the derivative
lands on $\beta(|\lambda|^{\delta_1} D(x,y))$ the argument we just used for the 
$|D(x,y)|^{\delta z}$ case works. One thing worth pointing out is that in these cases the presence of the $z$ in the 
exponent leads to an additional factor of $C|Im(z)|$ upon differentiation; however, the presence of the 
$e^{z^2}$ in the damping factor is more than enough to compensate. 

Lastly, we consider the case where the derivative lands on the factor $|H(x,y - \psi(x))|^{z({1 \over 2} - 
{1 \over d^*})}$. Since this factor was defined differently on the different $D_i$, $E_{ij}$, etc, we split
the square $S$ into its intersections with the $D_i''$, $E_{ij}''$, $F_{ij}''$, and $G_{ij}''$. For anything
other than a $G_{ij}''$ with $k_{ij}
> 2$, the damping factor is a power of $F^*(x,y - \psi(x))$. The directions $u$ and $v$ were defined so that
$F^*(x,y - \psi(x))^2$ has some nonvanishing higher order derivative in the $u$ and $v$ directions, so one can argue 
as above, breaking up the one-dimensional integration in the $u$ or $v$ variables into boundedly many 
intervals on which $F^*(x,y - \psi(x))$ is monotone. 

On a $G_{ij}''$ with $k_{ij} > 2$, the damping factor was defined as 
 $x^{M_i - {a_i + M_ib_i \over d^*}} \partial_{yy} f(x,y)$. We can actually assume that $u$ 
and $v$ are such that the $(k_{ij} - 2)$th $u$ and $v$ derivatives of $x^{M_i - {a_i + M_ib_i \over d^*}} 
\partial_{yy} f(x,y)$ are nonvanishing. To see why, first note that $(2.19)$ gives that the $(k_{ij} - 2)$th 
$y$-derivative of $\partial_{yy} f(x,y)$ is bounded below by $C x^{a_i + M_i(b_i - k_{ij})}$. On the other 
hand, by Lemma 2.3, (remembering that $M_i$ is always at least 1 in generic adapted coordinates) on 
$G_{ij}''$ we have $|\partial^{\alpha}f(x,y)| \leq C_{\alpha}^{a_i + M_ib_i - M_i|\alpha|}$. Using these 
facts with the product rule , if $u$ and $v$ are close enough to the $y$ 
direction, the $(k_{ij} - 2)$th derivative in the $u$ or $v$ direction of $x^{M_i - {a_i + M_ib_i \over d^*}}
\partial_{yy} f(x,y)$ will also be nonvanishing. Hence one can argue as in the previous paragraph and get
the same upper bounds as before. We have now considered all possible places the derivative lands,
concluding the proof of the desired upper bounds for $|J_1(\lambda)|$.

The bounds for $|J_2(\lambda)|$ are proven exactly as they were for $|J_1(\lambda)|$, replacing the roles of the $u$ and $v$ 
variables. The presence of the added condition $|\partial_uf(x,y) - a_1| \leq |\lambda|^{-{1 \over 3}}$ in 
the domain, which does not have an analogue above, does not interfere with any of the above estimates;
the condition that $\partial_u\partial_v^{k-1}D(x,y)$ is nonvanishing ensures that 
in any of the situations where one takes absolute values and does a Van der Corput type argument in the
$v$ direction, one still has boundedly many intervals. 

We now move on to $J_3(\lambda)$. Consider the level sets of $\partial_u f(x,y)$ and $\partial_v f(x,y)$. The 
gradients of 
both functions are bounded below in absolute value by $C|D(x,y)|$, which is at least ${1 \over 3}|\lambda|
^{-\delta_1}$ on the square $S$ if we chose the constant $c$ in the diameter $c|\lambda|^{-\delta_1}$ 
of the squares sufficiently small. As a result, if $c$ is small enough the level sets of both 
$\partial_u f(x,y)$ and $\partial_v f(x,y)$ do not self-intersect on $S$. Hence we may use $\partial_u f(x,y)$
and $\partial_v f(x,y)$ as coordinates on $S$. In particular, we may evaluate the measure of the set $S_3$
of $(4.5c)$ by changing into these coordinates in the integral of its characteristic function. The result is
$$|S_3| < C \min_S |D(x,y)|^{-1} |\lambda|^{-{2 \over 3}} \eqno (4.6)$$
So we conclude that $|S_3| < C'|\lambda|^{-{2 \over 3} + \delta_1}$. Since the integrand of $J_3(\lambda)$ is uniformly
bounded on $Re(z) = s$ for any $s > 1$, we conclude that 
$$|J_3(\lambda)| \leq C''|\lambda|^{-{2 \over 3} + \delta_1} < C''|\lambda|^{-{3 \over 5}} \eqno (4.7)$$
This gives the needed estimate. Adding the contributions from $|J_1(\lambda)|$, $|J_2(\lambda)|$, and $|J_3(\lambda)|$, we conclude
that the contribution to $|I_3(\lambda)|$ from the square $S$ is at most $C'''|\lambda|^{-{3 \over 5}}$,
and since $-{3 \over 5} < -{1 \over 2}$ we conclude that $|I_3(\lambda)|$ satisfies the bounds we need so 
long as $\delta_1$ was chosen sufficiently small.

\noindent {\medheading Estimating $|I_2(\lambda)|$}.

We focus our attention on the main term $I_2(\lambda)$, given by $(4.4b)$. We divide the domain of
$(4.4b)$ into squares of diameter $c|\lambda|^{-\delta_2}$, where $c$ and $\delta_2$ are small constants. 
For a given such square $S$,
denote the corresponding term of $I_2(\lambda)$ by $I_2^S(\lambda)$. We will show that if $c$ and $\delta_2$ are 
sufficiently small, then for any such $S$ we have $|I_2^S(\lambda)| < C|\lambda|^{-{1 \over 2} - \epsilon}$,
where $\epsilon$ is independent of $c$ and $\delta_2$, and $C$ is independent of $Im(z)$ for $Re(z) = s > 1$. 
Since there are at most $c'|\lambda|^{2\delta_2}$ squares, as long as we make sure $\delta_2 < {\epsilon 
\over 2}$, this is enough to show that $I_2(\lambda)$ itself satisfies the 
bounds needed for Theorem 4.1. This subdivision into squares is useful because it allows us to 
replace $D(x,y)$ by a polynomial approximation of bounded degree which is therefore 
piecewise monotone in a direction in which we are integrating by parts, enabling us to use Van der
Corput type arguments in such a direction.

We now perform this polynomial replacement. For a given $S$ and positive integer $N$, let $D_S^N(x,y)$ be 
the polynomial in $x$ and $y$ consisting of the sum of the terms of degree at most $N$ of 
$D(x,y)$'s Taylor expansion centered about the center of $S$. Thus on $S$ we have
$$|D_S^N(x,y) - D(x,y)| < C|\lambda|^{-\delta_2 N} \eqno (4.8)$$
As a result, on $S$ we have
$$|\big( \alpha(|\lambda|^{\delta_1} D(x,y)) - \alpha(|\lambda|^{N_1} D(x, y)\big)
- \big( \alpha(|\lambda|^{\delta_1} D_S^N(x,y)) - \alpha(|\lambda|^{N_1} D_S^N(x, y))\big)| < 
C|\lambda|^{N_1 - \delta_2 N}$$
In particular, if $N$ is chosen large enough we can make the exponent $N_1 - \delta_2 N$ appearing in $(4.9)$
less than -1. Consequently, for the purposes of the analysis of $I_2^S(\lambda)$ we may replace $D(x,y - 
\psi(x))$ by $D_S^N(x,y)$ in the $\big( \alpha(|\lambda|^{\delta_1} D(x,y) - 
\alpha(|\lambda|^{N_1} D(x, y))\big)$ factor; the difference will contribute no more than $C
|\lambda|^{-1}$ to $I_2^S(\lambda)$, and adding over all squares gives a result smaller than the bounds
needed in Theorem 4.1.

We can do something similar for the $|D(x,y)|^{\delta z}$ factor. Namely, suppose $N$ is taken 
large enough that in $(4.8)$ we have 
$$|D_S^N(x,y) - D(x,y)| < C|\lambda|^{-2N_1}$$
Then since $|D(x, y)| \geq {1 \over 2}|\lambda|^{-N_1}$ when the integrand of $(4.4b)$ is nonzero,
if $|\lambda|$ is large enough we may use the Taylor expansion of $|x|^{\delta z}$ about $x = 
D(x, y)$ to obtain  
$$\big||D(x,y)|^{\delta z} - |D_S^N(x,y)|^{\delta z}\big| < C|Im(z)| 
|\lambda|^{-N_1(\delta Re(z) - 1) - 2N_1} \eqno (4.9)$$
As a result, since $Re(z) > 1$, as long as $N_1 > 1$, we have an estimate
$$\big||D(x,y)|^{\delta z} - |D_S^N(x,y)|^{\delta z}\big| < C|Im(z)| 
|\lambda|^{-1} \eqno (4.10)$$
The $e^{z^2}$ is more than enough to take care of the $|Im(z)|$ factor in $(4.10)$, and the exponent 
$-1$ is less than $-{1 \over 2}$. Consequently, we may replace $|D(x,y)|^{\delta z}$
by $|D_S^N(x,y)|^{\delta z}$ in the analysis of $I_2^S(\lambda)$; the difference added over all squares $S$
contributes less than the bounds needed for Theorem 4.1. 

We have now shown that for the purposes of our
future arguments, we may adjust our notation and assume $I_2^S(\lambda)$ is given by
$$I_2^S(\lambda)= e^{z^2} \int e^{-i\lambda_1 f(x,y) - i\lambda_2x - i\lambda_3 y}|H(x,y - \psi(x))|^z
|D_S^N(x,y)|^{\delta z}$$
$$\times \big( \alpha(|\lambda|^{\delta_1} D_S^N(x,y)) - \alpha(|\lambda|^{N_1} D_S^N(x,y))\big)
\phi^*(x,y) \,dx\,dy \eqno (4.11)$$

We divide the domain of integration of $(4.11)$ into the intersections of $S$ with the $D_i'', E_{ij}'',
F_{ij}'',$ and $G_{ij}''$  and denote the corresponding term of $I_2(\lambda)$
by $I_2^{D_i}(\lambda), I_2^{E_{ij}}(\lambda)$, $I_2^{F_{ij}}(\lambda)$, and $I_2^{G_{ij}}(\lambda)$. (Recall
$D_i'' = \{(x,y): (x, y - \psi(x)) \in D_i\}$ with analogous definitions for the other regions). We 
suppress the $S$ since the bounds we will prove, given in the statement of Theorem 4.1, are independent of
$S$. We will only consider those regions for which $x > 0$ as the $x < 0$ ones are entirely analogous.
We now focus our attention on the analysis of the $I_2^{D_i}(\lambda)$.

\noindent {\bf Bounds for $|I_2^{D_i}(\lambda)|$}.

\noindent Recalling that $|H(x,y)| = F^*(x,y)^{{1 \over 2} - {1 \over d^*}}$ on a $D_i$, if we
change coordinates from $(x,y)$ to $(x,y + \psi(x))$ in $(4.11)$ we obtain
$$I_2^{D_i}(\lambda) =  e^{z^2} \int_{S \cap D_i} e^{-i\lambda_1 F(x,y) - i\lambda_2x - i\lambda_3 
(y + \psi(x))}F^*(x,y)^{({1 \over 2} - {1 \over d^*})z} |D_S^N(x,y + \psi(x))|^{\delta z}$$
$$\big( \alpha(|\lambda|^{\delta_1} D_S^N(x,y + \psi(x))) - \alpha(|\lambda|^{N_1} D_S^N(x, y + \psi(x)))\big)
\phi^{**}(x,y) \,dx\,dy \eqno (4.12)$$
Here $\phi^{**}(x,y)$ denotes a new cutoff function on a neighborhood of the origin, and $F(x,y)
= f(x, y + \psi(x))$ is in generic adapted coordinates not satisfying the exceptional situations of Theorem
1.1. We slightly abuse 
notation in $(4.12)$ in that $S$ now denotes the square in the new coordinates. We now decompose the 
domain of $(4.12)$ into dyadic rectangles. We only consider those rectangles in the upper right quadrant
as the other quadrants are done the same way. For a given dyadic rectangle $J_{kl} = [2^{-k-1},2^{-k}] 
\times [2^{-l-1},2^{-l}]$, we use the shorthand by $I_{kl}$ to denote the corresponding term of 
$I_2^{D_i}(\lambda)$, given by
$$I_{kl} =  e^{z^2} \int_{S \cap D_i \cap J_{kl}} e^{-i\lambda_1 F(x,y) - i\lambda_2x - i\lambda_3 
(y + \psi(x))}F^*(x,y)^{({1 \over 2} - {1 \over d^*})z}|D_S^N(x,y + \psi(x))|^{\delta z}$$
$$\big( \alpha(|\lambda|^{\delta_1} D_S^N(x,y + \psi(x))) - \alpha(|\lambda|^{N_1} D_S^N(x, y + \psi(x)))\big)
\phi^{**}(x,y) \,dx\,dy \eqno (4.13)$$

We will analyze $(4.13)$ by imitating the proof of Van der Corput's lemma in the $y$ direction. Our 
objective is to show that $(4.13)$ is bounded by $C(1 +|\lambda|)^{-{1 \over 2} - \epsilon}$ as in the 
statement of Theorem 4.1. The second $y$ derivative of the phase function
in $(4.13)$ is given by $\lambda_1 \partial_{yy} F(x,y)$, 
and by $(2.14)$, if the vertex of $N(F)$ corresponding to $D_i$ is written as $(a_i,b_i)$, then on $D_i$ we 
have $|\partial_{yy} F(x,y)| > c|x|^{a_i}|y|^{b_i - 2}$. Since $x \sim 2^{-k}$ and $y \sim 2^{-l}$ on 
$J_{kl}$ we can write this as 
$$|\partial_{yy} F(x,y)| > c'{1 \over (2^{-l})^2} (2^{-k})^{a_i} (2^{-l})^{b_i} \eqno (4.14)$$
As in the proof of the Van der Corput theorem for
functions with nonvanishing second derivative, we will split the integral $(4.13)$ into two parts. The
first is the part where $|\lambda_1\partial_y F(x,y) + \lambda_3| < |\lambda|^{{1 \over 2}} 2^{{-a_i k  - (b_i - 2) l \over 2}}$,
and the second is the part where $|\lambda_1\partial_y F(x,y) + \lambda_3| \geq |\lambda|^{{1 \over 2}} 
2^{{-a_i k  - (b_i - 2)l \over 2}}$. Call the resulting integrals $K_1$ and $K_2$, so that $K_1 + K_2 = 
I_{kl}$. We will bound $K_1$ by taking absolute values and integrating, and $K_2$ by performing an 
integration by parts.

We start with $K_1$. The integrand of $(4.13)$ is bounded in absolute value by a constant times 
$F^*(x,y)^{Re(z)({1 \over 2} - {1 \over d^*})}|D_S^N(x,y + \psi(x))|^{\delta Re(z)}$. By $(2.20)$, 
$F^*(x,y) < C|x^{a_i}y^{b_i}| \leq C2^{-ka_i -lb_i}$, and on the 
domain of $(4.13)$ we have $|D_S^N(x,y + \psi(x))|^{\delta z} < C'|\lambda|^{-\delta \delta_1 Re(z)}$
Hence if $s$ denotes $Re(z)$, the integrand of $(4.13)$ is at most 
$$C''|\lambda|^{-\delta\delta_1 s}2^{(-ka_i -lb_i)s({1 \over 2} - {1 \over d^*})} \eqno (4.15)$$
Since by $(4.14)$ the absolute value of the $y$-derivative of $\lambda_1\partial_y F(x,y) + \lambda_3$ is 
at least $c'|\lambda| 2^{-k a_i - l(b_i - 2)}$ we have
$$|\{y: |\lambda_1\partial_y F(x,y) + \lambda_3| < |\lambda|^{{1 \over 2}}2^{{-a_i k  - (b_i - 2) 
l \over 2}}\}| < |\lambda|^{-{1 \over 2}} 2^{{k a_i   + l(b_i - 2) \over 2}} \eqno (4.16)$$
Thus bounding the $y$ integral of $K_1$ by $(4.15)$ times the measure $(4.16)$ and then integrating the
result in $x$, we obtain
$$|K_1| < C'''|\lambda|^{-{1 \over 2} - \delta\delta_1 s} 2^{(-k a_i  - lb_i)({s - 1 \over 2} - 
{s \over d^*})}2^{-k - l} \eqno (4.17)$$
We now turn to $K_2$ and show that $K_2$ also satisfies the upper bounds of $(4.17)$. We integrate the
integrand in $(4.13)$ by parts
in $y$, integrating the factor $(\lambda_1\partial_y F(x,y) + \lambda_3)e^{-i\lambda_1 F(x,y) - i\lambda_2x - 
i\lambda_3 (y - \psi(x))}$ and differentiating ${1 \over \lambda_1\partial_y F(x,y) + \lambda_3}$ times
the rest of the integrand. We get several terms depending on where the derivative lands. If the derivative
lands on $\phi^*(x,y)$, the absolute value of the integrand in the resulting term is bounded by 
$$CF^*(x,y)^{s({1 \over 2} - {1 \over d^*})} |D_S^N(x,y + \psi(x))|^{\delta s} |\lambda|^{-{1 \over 2}} 
2^{{a_i k  + (b_i - 2)l \over 2}}\eqno (4.18)$$
Bounding $F^*(x,y) < C2^{-ka_i - lb_i}$ and  $|D_S^N(x,y + \psi(x))| < C|\lambda|^{-\delta_1}$ as 
in the analysis of $K_1$, we get that $(4.18)$ is bounded by
$$C'2^{(-ka_i - lb_i)({s-1 \over 2} - {s \over d^*})}2^{-l}|\lambda|^{-{1 \over 2} - \delta \delta_1 s}
\eqno (4.19)$$
Integrating $(4.19)$ over $S \subset I_{kl}$ multiplies this by at most $C2^{-k - l}$, so the resulting term
is at most
$$C'2^{(-ka_i - lb_i)({s-1 \over 2} - {s \over d^*})}2^{-k - 2l}|\lambda|^{-{1 \over 2} - 
\delta \delta_1 s} \eqno (4.20)$$
Note this is better than the estimate $(4.17)$. We next consider the case where the $y$-derivative lands on
the ${1 \over \lambda_1\partial_y F(x,y) + \lambda_3}$ factor, turning it into $-{\lambda_1\partial_{yy} F(x,y)
\over (\lambda_1\partial_y F(x,y) + \lambda_3)^2}$. We take absolute values and integrate in the $y$ 
variable as in the proof of the Van der Corput lemma, bounding the other factors as was done for
$(4.19)$. Since by $(4.14)$ the function $\partial_{yy} F(x,y)$ is never zero on the domain of integration,
we have at most finitely many intervals of integration on each of which 
$|{\lambda_1\partial_{yy} F(x,y) \over (\lambda_1\partial_y F(x,y) + \lambda_3)^2}|$ integrates back into
$\pm{1 \over \lambda_1\partial_y F(x,y) + \lambda_3}$. Hence the resulting term, as well as the endpoint terms, will be bounded by $(4.20)$,
except divided by the $y$-width $2^{-l}$. We conclude that this term is bounded by $(4.17)$, namely
$$C''2^{(-ka_i - lb_i)({s-1 \over 2} - {s \over d^*})}2^{-k - l}|\lambda|^{-{1 \over 2} - 
\delta \delta_1 s} \eqno (4.21)$$
If the $y$-derivative lands on either the $|D_S^N(x,y + \psi(x))|^{\delta z}$ or 
$\big( \alpha(|\lambda|^{\delta_1} D_S^N(x,y + \psi(x))) - \alpha(|\lambda|^{N_1} D_S^N(x,y + \psi(x)))\big)$ 
factors one estimates the resulting term in very much the same way; the fact that $D_S^N(x,y)$ is a 
polynomial and $\alpha$ is monotone ensures that the Van der Corput lemma proof still applies and we will 
have boundedly many intervals of
integration on which the appropriate derivative is nonvanishing. Similarly, since $F^*(x,y)^2$ is a polynomial,
one can deal with the term where the derivative lands on the damping factor $F^*(x,y)^{z({1 \over 2} -
{1 \over d^*})}$ in a similar fashion. It should be pointed out that in taking these derivatives we do incur
a factor of $C|Im(z)|$, but this is more than compensated for by the $e^{z^2}$ factor.
Hence we once again get the upper bound $(4.21)$. Adding all terms
together, we see that $|K_2|$ and therefore $|I_{kl}|$ is bounded by $(4.17)$, the estimate we need.

We rewrite $(4.17)$ in an especially useful form. Recall that by $(2.20)$, on $S$ we have
$C2^{(-ka_i - lb_i)} < F^*(x,y) < C'2^{(-ka_i - lb_i)}$. So we have just shown that 
$$|I_{kl}| < C\int_{J_{kl}} |\lambda|^{-{1 \over 2} - 
{\delta \delta_1 s \over 2}}F^*(x,y)^{{s-1 \over 2} - {s \over d^*}}\,dx\,dy \eqno (4.22)$$
We now break into cases $d^* \leq 2$, and $d^* > 2$, starting with the latter. Adding $(4.22)$ 
over all rectangles, we obtain that $|I_2^{D_i}(\lambda)|$ is at most
$$C |\lambda|^{-{1 \over 2} - {\delta \delta_1 s \over 2}} \int_{[0,1] \times [0,1]}F^*(x,y)^{{s-1 
\over 2} - {s \over d}}\,dx\,dy \eqno (4.23)$$
Note that if $s > 1$, then ${s-1 \over 2} - {s \over d} > -{1 \over d}$, and thus since $F^*(x,y)^t$ is 
integrable over $[0,1] \times [0,1]$ for all $t > -{1 \over d}$, the integral in
$(4.23)$ is finite and we obtain that $|I_2^{D_i}(\lambda)|$
is bounded by $C |\lambda|^{-{1 \over 2} - {\delta \delta_1 s \over 2}}$. Since the exponent here is less 
than $-{1 \over 2}$, this gives what is needed for Theorem 4.1.

\noindent Moving on to the $d^* = 2$ case, $(4.22)$ becomes
$$|I_{kl}| < C\int_{J_{kl}} |\lambda|^{-{1 \over 2} - 
{\delta \delta_1 s \over 2}}F^*(x,y)^{-{1 \over 2}}\,dx\,dy \eqno (4.24)$$
Since the damping factor is just $|D_S^N(x,y + \psi(x))|^{\delta z}$ when $d^* = 2$, from $(4.13)$ we get
$$I_{kl} =  e^{z^2} \int_{S \cap D_i \cap J_{kl}} e^{-i\lambda_1 F(x,y) - i\lambda_2x - i\lambda_3 
(y + \psi(x))}|D_S^N(x,y + \psi(x))|^{\delta z}$$
$$\times \big( \alpha(|\lambda|^{\delta_1} D_S^N(x,y + \psi(x))) - \alpha(|\lambda|^{N_1} D_S^N(x, y + \psi(x)))\big)
\phi^{**}(x,y) \,dx\,dy \eqno (4.25)$$
Note that due to the cutoff and the presence of the $|D_S^N(x,y + \psi(x))|^{\delta z}$ in the integrand of
$(4.25)$, this integrand is at most $|\lambda|^{-\delta \delta_1s}$. So just by taking absolute
values and integrating we get
$$|I_{kl}| < C|\lambda|^{-\delta \delta_1 s}2^{-k - l} \eqno (4.26a)$$
$$< C|\lambda|^{-\delta \delta_1 s \over 2}2^{-k - l} \eqno (4.26b)$$
Combining this with $(4.24)$, we get
$$|I_{kl}| < C |\lambda|^{-\delta \delta_1 s \over 2} \int_{[2^{-k-1}, 2^{-k}] \times 
[2^{-l-1}, 2^{-l}]}\min (1, |\lambda F^*(x,y)|^{-{1 \over 2}})\,dx\,dy \eqno (4.27)$$
Adding this up over all $j$ and $k$, we obtain that $|I_2^{D_i}(\lambda)|$ is at most 
$$C |\lambda|^{-\delta \delta_1 s \over 2} \int_{[0,1] \times [0,1]} \min (1, |\lambda F^*(x,y)|
^{-{1 \over 2}})\,dx\,dy\eqno (4.28)$$
Since $(d,d) \in N(F)$, $(d,d)$ is a convex combination of vertices of $N(F)$. So since $F^*(x,y)$ is 
comparable to the 
sum of $|x^ay^b|$ over vertices $(a,b)$ of $N(F)$, we have $F^*(x,y) > C|x^dy^d|$. Since
we are assuming $d \leq 2$ here, we conclude that $F^*(x,y) > Cx^2y^2$ and as a result $(4.28)$ is bounded by 
$$C |\lambda|^{-\delta \delta_1 s \over 2} \int_{[0,1] \times [0,1]} \min \bigg(1, { C' \over |\lambda|^
{1 \over 2}xy}\bigg)\,dx\,dy\eqno (4.29)$$
A direct calculation reveals that the right hand side is bounded above by $C'' |\lambda|^{-{1 \over 2}}
(\ln |\lambda|)^2$ (The integral over $[|\lambda|^{-{1 \over 2}}, 1] \times [|\lambda|^{-{1 \over 2}},1]$ is bounded by 
a constant times the integral of ${ 1 \over |\lambda|^{1 \over 2}xy}$ over this region, while the integral 
over the remaining region is bounded by its area). As a result, $(4.29)$ is bounded by
$$C''|\lambda|^{-{1 \over 2} -{\delta \delta_1 s \over 2}}(\ln|\lambda|)^2 \eqno (4.30)$$
Since the exponent here is less than $-{1 \over 2}$ we have
proved the desired bounds for the $|I_2^{D_i}(\lambda)|$.

\noindent {\bf Bounds for $|I_2^{E_{ij}}(\lambda)|$}.

\noindent Note that $I_2^{E_{ij}}(\lambda)$ is given by  
$$I_2^{E_{ij}}(\lambda)= e^{z^2} \int_{S \cap E_{ij}''} e^{-i\lambda_1 f(x,y) - i\lambda_2x - i\lambda_3 y}
|F^*(x,y - \psi(x))|^{z({1 \over 2} - {1 \over d^*})} |D_S^N(x,y)|^{\delta z}$$
$$\big(\alpha(|\lambda|^{\delta_1} D_S^N(x,y)) - \alpha(|\lambda|^{N_1} D_S^N(x,y))\big)
\phi^*(x,y) \,dx\,dy \eqno (4.31)$$ 
As we did with $I_2^{D_i}(\lambda)$, we break the domain of integral $(4.31)$ into  
rectangles $J_{kl} =  [2^{-k-1}, 2^{-k}] \times [2^{-l-1}, 2^{-l}]$. Denote the corresponding term of
$(4.31)$ by $I_{kl}$, so that $\sum_{kl}I_{kl} = I_2^{E_{ij}}(\lambda)$. 

\noindent Note that by $(2.21)$ we have
$$F^*(x,y) > Cx^{a_i + M_ib_i} \eqno (4.32)$$
As before $(a_i,b_i)$ denotes the upper vertex of $e_i$. Recall that $E_{ij}$ lies between $y = 
(r - \eta)x^{M_i}$ and $y = (r + \eta)x^{M_i}$ for some $r$ and $\eta$ such that $F_{e_i}(1,r)
\neq 0$, and that by definition of $E_{ij}$, $\psi(x)$ has a zero of order at least $M_i$ at $x = 0$.
Consequently, $|y - \psi(x)| < Cx^{M_i}$ on $E_{ij}$. Thus by $(2.11)$, $F^*(x,y - \psi(x)) < Cx^{a_i +
b_iM_i}$ on $E_{ij}$. Combining with $(4.32)$ we get
$$ F^*(x,y - \psi(x)) < C F^*(x,y) \eqno (4.33a)$$
By Lemma 2.4, on the domain of $(4.31)$ we have
$$|\partial_{xx} f(x,y)| > C{1 \over x^2} x^{a_i}(x^{M_i})^{b_i} \eqno (4.33b)$$
Equation $(4.33a)$ shows that the damping function $F^*(x,y - \psi(x))$ satisfies the same upper bounds
that the damping function $F^*(x,y)$ did in the $D_i$ case.
Equation $(4.33b)$ shows the same thing for the phase (cf $(4.14)$), reversing the roles of the $x$ and $y$
derivatives. Furthermore, the functions that need to be piecewise monotone in $x$ with boundedly many pieces
in order to perform the Van der Corput argument do satisfy this; $D_S^N(x,y)$ is a polynomial and the second
$x$ derivative of $-i\lambda_1 f(x,y) - i\lambda_2x - i\lambda_3 y$ is nonvanishing by $(4.33b)$. 
Hence by repeating the $D_i$ argument, reversing the roles of the $x$ and $y$ variables, we get that
$I_{kl}$ is bounded by $(4.17)$. Adding this up like before gives that
as in $(4.30)$, $|I_2^{E_{ij}}(\lambda)|$ is 
bounded by $C''|\lambda|^{-{1 \over 2} -{\delta \delta_1 s \over 2}}(\ln|\lambda|)^2$, the estimate we need.

\noindent {\bf Bounds for $|I_2^{F_{ij}}(\lambda)|$}.

Recall the set $F_{ij}$ is of the form
$\{(x,y): 0 < x < \eta, |y - rx^{M_i}| < \nu|x|^{M_i} \}$, where $F_{e_i}(1,y)$ has a zero of order 1 at
$y = r$. Define $G(x,y) = F(x,y + rx^{M_i})$. Thus $G(x,y)$ is a function on the set $H_{ij} = 
\{(x,y): 0 < x < \eta, |y| < \nu|x|^{M_i} \}$ such that $G_{e_i}(1,y)$ has a zero of order 1 at $y = 0$. Thus
$N(\partial_y G)$ has an edge with equation $x + M_i y = a_i + M_ib_i - M_i$ that intersects the $x$ axis.
Conseqently, $N({\partial^2 G \over \partial x \partial y})$ has an edge with equation $x + M_i y = a_i +
M_ib_i - M_i - 1$ intersecting the $x$ axis. Hence assuming $\eta$ was chosen sufficiently small, by 
Lemma 2.3 we may conclude that on $H_{ij}$ we have
$$\bigg|{\partial^2 G \over \partial x \partial y}(x,y)\bigg| > C x^{a_i + M_ib_i - M_i - 1}$$
We rewrite this as 
$$ \bigg|{\partial^2 G \over \partial x \partial y}(x,y)\bigg| > C {1 \over x (x^{M_i})}x^{a_i}(x^{M_i})^{b_i} \eqno (4.34)$$
Letting $\tilde{\psi}(x) = \psi(x) + rx^{M_i}$, we do a change of variables from $y$ to $y + 
\tilde{\psi}(x)$ and write $I_2^{F_{ij}}(\lambda)$ as 
$$I_2^{F_{ij}}(\lambda)= e^{z^2} \int_{S \cap H_{ij}} e^{-i\lambda_1 G(x,y) - i\lambda_2x - i\lambda_3 
(y + \tilde{\psi}(x))}|F^*(x,y + rx^{M_i})|^{z({1 \over 2} - {1 \over d^*})}$$
$$\times |D_S^N(x,y + \tilde{\psi}(x))|^{\delta z} \big(\alpha(|\lambda|^{\delta_1} D_S^N(x,y+  
\tilde{\psi}(x))) - \alpha(|\lambda|^{N_1} D_S^N(x,y + \tilde{\psi}(x)))\big)\phi^{***}(x,y) \,dx\,dy
\eqno (4.35)$$ 
As with the $D_i$, the $S$ under the integral symbol now denotes the square in the new coordinates.
By $(2.21)$, on $H_{ij}$ we have
$$F^*(x,y + rx^{M_i}) < C x^{a_i}(x^{M_i})^{b_i} \eqno (4.36)$$
We now break the domain of integration of $(4.35)$ up into 
rectangles $J_k$ of the form $[2^{-k - 1}, 2^{-k}] \times [-\nu 2^{-kM_i}, \nu 2^{-kM_i}]$, and let
$I_k'$ the the portion of $(4.35)$ coming from $J_k$.
Equation $(4.36)$ shows that the damping function $F^*(x,y + rx^{M_i})$ in $(4.35)$ satisfies
the same upper bounds the damping function did on the $I_{kM_k}$ rectangle for the the $D_i$. (The 
$x \sim 2^{-k}$ rectangle of the "lower edge" of $D_i$). As for the 
phase, instead of having a lower bound on a second $y$ derivative as in $(4.14)$, we have the substitute $(4.34)$. We still
may argue as for the $I_{kM_k}$ rectangle in the $D_i$ case, but with one difference. In the analysis of 
the term called $K_1$ below $(4.14)$, instead of bounding the
measure of a sublevel set of $|\lambda_1\partial_y G(x,y) + \lambda_3|$ in the $y$-variable and integrating 
with respect to $x$, one bounds the measure of the same sublevel set in the $x$ variable using $(4.34)$ and
then integrates the result with respect to $y$.

Furthermore, all relevant factors are piecewise monotone with boundedly many pieces. The function 
$D_S^N(x,y + \tilde{\psi}(x))$ is a polynomial in $y$ of bounded degree, as is $F^*(x,y + rx^{M_i})^2$, while since $G(x,y)$ is just a
$y$ shift of $F(x,y)$ by $\psi(x)$, if $(0,a)$ denotes the upper vertex of $N(F)$ then $\partial_y G(x,y)$
has nonvanishing $(a - 1)$th $y$ derivative. 

Hence after making the above adjustment to the $I_{kM_ik}$ argument of the $D_i$ case, for a given $k$ we get the bounds 
$(4.17)$ for $I_k'$. (The arguments there did not require $l$ to be an integer). Adding over all $k$,
as for the $|I_2^{D_i}(\lambda)|$ we get that 
$|I_2^{F_{ij}}(\lambda)|$ is bounded by $C''|\lambda|^{-{1 \over 2} -{\delta \delta_1 s \over 2}}
(\ln|\lambda|)^2$, the needed estimate.

\noindent {\bf Bounds for $|I_2^{G_{ij}}(\lambda)|$}.

For the $I_2^{G_{ij}}(\lambda)$, we separate the $k_{ij} = 2$ and $k_{ij} > 2$ cases as the
damping factors are different in these two situations. First, we suppose $k_{ij} = 2$. Then 
$I_2^{G_{ij}}(\lambda)$ is given by
$$I_2^{G_{ij}}(\lambda)= e^{z^2} \int_{S \cap G_{ij}} e^{-i\lambda_1 F(x,y) - i\lambda_2x - i\lambda_3 
(y + \psi(x))}F^*(x,y)^{z({1 \over 2} - {1 \over d^*})}|D_S^N(x,y + \psi(x))|^{\delta z}$$
$$\big(\alpha(|\lambda|^{\delta_1} D_S^N(x,y +  \psi(x))) - 
\alpha(|\lambda|^{N_1} D_S^N(x,y + \psi(x)))\big)\phi^{**}(x,y) \,dx\,dy \eqno (4.37)$$ 
Observing that $|y| < Cx^{M_i}$ on $G_{ij}$, we divide the domain of $(4.37)$ into rectangles $J_k$ of the
form $[2^{-k-1},2^{-k}] \times [C_0 2^{-kM_i}, C_12^{-kM_i}]$, and let $I_k'$ be the corresponding piece of
$I_2^{G_{ij}}(\lambda)$, so that $\sum_k I_k' = I_2^{G_{ij}}(\lambda)$.

Note that the integrand in $(4.37)$ is the same as that of $(4.13)$ for the $D_i$ case. In particular, 
the damping function is the same as in the $D_i$ case. Also, by $(2.19)$ on $G_{ij}$ we have the following
analogue of $(4.14)$:
$$|\partial_{yy} F(x,y)| > C{1 \over (x^M_i)^2}x^{a_i}(x^{M_i})^{b_i} \eqno (4.38)$$
As a result, all estimates used in the $D_i$ case for the $I_{kl}$ rectangle, setting $l = kM_i$ (the 
lower edge of $D_i$) hold for 
the term $I_k'$. Thus $|I_k'|$ is bounded by $C|I_{kM_ik}|$, and adding over all $k$ we recover
$C|\lambda|^{-{1 \over 2} -{\delta \delta_1 s \over 2}}(\ln|\lambda|)^2$ as an upper bound for 
$|I_2^{G_{ij}}(\lambda)|$. This completes the proof for the $k_{ij} = 2$ case.

\noindent We may now assume $k_{ij} > 2$, focusing our attention for now on the case when $N(F)$ has 
multiple vertices. Here, $I_2^{G_{ij}}(\lambda)$ is given by
$$I_2^{G_{ij}}(\lambda)= e^{z^2} \int_{S \cap G_{ij}} e^{-i\lambda_1 F(x,y) - i\lambda_2x - i\lambda_3 
(y + \psi(x))}\big[x^{M_i - {a_i + M_ib_i \over d}}|{\partial^2F \over \partial y^2}(x,y)|^{1 \over 2}\big]
^z$$
$$|D_S^N(x,y + \psi(x))|^{\delta z}\big(\alpha(|\lambda|^{\delta_1} D_S^N(x,y +  \psi(x))) - 
\alpha(|\lambda|^{N_1} D_S^N(x,y + \psi(x)))\big)\phi^{**}(x,y) \,dx\,dy \eqno (4.39)$$
We divide the domain of $(4.39)$ into rectangles $J_k$ as in the above $k_{ij} = 2$ case.
and again let $I_k'$ be the corresponding piece of $I_2^{G_{ij}}(\lambda)$. Observe that by $(2.23)$, there
is some $C_0$ such that the magnitude of the bracketed expression in 
$(4.39)$ (which is the same as the $H(x,y)$ in $(2.23)$) is bounded by $C_0(x^{a_i + M_ib_i})^{{1 \over 2} - 
{1 \over d}}$. Thus we may write $I_k' = \sum_{l = 0}^{\infty} P_{kl}$, where $P_{kl}$ is the portion of the 
integral over $J_k$ where $|H(x,y)|$
is between $2^{-l+1}C_0(x^{a_i + M_ib_i})^{{1 \over 2} - {1 \over d}}$ and $2^{-l}C_0(x^{a_i + M_ib_i})
^{{1 \over 2} - {1 \over d}}$. We will now bound each $P_{kl}$. To this end, note
that on the domain of $P_{kl}$, by the definition of $H(x,y)$ and the $P_{kl}$ we have
$$C_0^22^{-2l - 2}(x^{a_i + M_ib_i})^{1 - {2\over d}} < x^{2M_i - {2a_i + 2M_ib_i \over d}}
|{\partial^2F \over \partial y^2}(x,y)| < C_0^22^{-2l}(x^{a_i + M_ib_i})^{1 - {2\over d}} \eqno (4.40)$$
Solving for $|{\partial^2F \over \partial y^2}(x,y)|$, we get
$$C_12^{-2l}{1 \over (x^{M_i})^2}x^{a_i + M_ib_i} < |{\partial^2F \over \partial y^2}(x,y)|
< C_1'2^{-2l}{1 \over (x^{M_i})^2}x^{a_i + M_ib_i} \eqno (4.41)$$
One now bounds $P_{kl}$ by integrating by parts in $y$ in the portion of $(4.39)$ corresponding to $P_{kl}$. One 
proceeds exactly as for the $I_{kM_ik}$ term of the $D_i$ (the $x \sim 2^{-k}$ rectangle of the "lower
edge" of $D_i$), except instead of using $|{\partial^2F \over 
\partial y^2} (x,y)| > C {1 \over (x^{M_i})^2} x^{a_i+ M_ib_i}$ from $(4.14)$ one uses $(4.41)$. This gives us an
additional factor of $C2^l$ in the resulting bounds for the integral. This however is compensated by the damping factor, which
by the definition of $P_{kl}$ is bounded by $C2^{-l Re(z)}$ times the damping factor used for the $I_{kM_ik}$
term in the $D_i$ case.
Thus the overall integral is bounded by $C2^{l(1 - Re(z))}$ times what is obtained for the $I_{kM_ik}$ term in 
the $D_i$ case. We do not have to worry about whether each factor in $(4.37)$ is boundedly piecewise 
monotone in $y$ in our integrations by parts; the
only new element in this regard is ${\partial^2F \over \partial y^2}(x,y)$, whose $(k_{ij} - 2)$th $y$ 
derivative is nonvanishing. 

Since $Re(z) > 1$, we conclude $\sum_l P_{kl}$ is bounded by a constant times the estimate
obtained for the $I_{kM_ik}$ term in the $D_i$ situation, and adding this over all $k$ gives
$$|I_2^{G_{ij}}(\lambda)| \leq \sum_{kl}|P_{kl}| < C''|\lambda|^{-{1 \over 2} -{\delta \delta_1 s \over 2}}
(\ln|\lambda|)^2 \eqno (4.42)$$
This is the estimate we seek. The above argument was for when $N(F)$ has multiple vertices, but when $N(F)$
just has one vertex the following simplified version of this argument works. In the one vertex situation, $|H(x,y)| = |{\partial^2F 
\over \partial y^2}(x,y)|$. This time we let $P_l$ be the portion of the integral defining $I_2^{G_{ij}}
(\lambda)$ over the set where $|H(x,y)|$ is between $C_02^{-l-1}$ and $C_02^{-l}$, where $C_0$ denotes the
maximum value of $|H(x,y)|$. Like above, for $P_l$ the decreased second $y$ derivative of the phase gives an
additional factor of $C2^l$ which is more than compensated by the additional
$C2^{-l Re(z)}$ factor coming from the damping function. Adding over all $l$, we recover $(4.42)$.
This completes the proof of the bounds for the $|I_2^{G_{ij}}(\lambda)|$, which in turn completes the proof
of Theorem 4.1.

\noindent {\bf The proof of Theorem 1.1.}

We may now finish the proof of Theorem 1.1 in short order. First suppose $d(F) > 2$. For any $\eta > 0$, 
Theorem 3.1 says that on the line $Re(z) = -{2 \over d(F) - 2} + \eta$, $M_z$ is bounded on $L^{\infty}$ with 
uniform constant, while Theorem 4.1 in conjunction with Theorem 1.2 says that on $Re(z) = 1 + \eta$, $M_z$ is bounded on $L^2$ with 
uniform constant. Using interpolation for maximal operators (see Ch. 11 of [St2]), we have that $M_0$ is 
bounded on $L^{d(F) + \eta'}$ where $\eta' \rightarrow 0$ as $\eta \rightarrow 0$. Thus we conclude $M_0$
is bounded on $L^p$ for all $p > d(F)$. Since $d(F) = h(q_0) = \max(2, h(q_0))$, this gives Theorem 1.1 for 
$d(F) > 2$.

On the other hand, if $d(F) = 2$, Theorem 3.1 says that on any vertical line $Re(z) = s$, $M_z$ is bounded 
on $L^{\infty}$ with uniform constant, and Theorem 4.1 still applies on a line $Re(z) = 1 + \eta$. Thus
interpolation now gives the result obtained by letting $d(F)$ approach 2 in the previous paragraph, namely
that $M$ is bounded on $L^p$ for $p > 2 = \max(h(q_0), 2)$. This completes the proof of Theorem 1.1.

\noindent {\bf 5. References.}

\noindent [B] J. Bourgain, {\it Averages in the plane over convex curves and maximal operators},
J. Anal. Math. {\bf 47} (1986), 69--85. \parskip = 4pt\baselineskip = 4pt

\noindent [CoMa] M. Cowling, G. Mauceri, {\it Inequalities for some maximal functions. II},  Trans. Amer. 
Math. Soc. {\bf 298} (1986),  no. 1, 341--365. 

\noindent [G1] M. Greenblatt, {\it Maximal averages over hypersurfaces and the Newton polyhedron}, submitted.

\noindent [Gr] A. Greenleaf, {\it Principal curvature and harmonic analysis}, 
Indiana Univ. Math. J. {\bf 30} (1981), no. 4, 519--537.

\noindent [IkKeMu1] I. Ikromov, M. Kempe, and D. M\"uller, {\it Damped oscillatory integrals and boundedness of
maximal operators associated to mixed homogeneous hypersurfaces} (English summary) Duke Math. J. {\bf 126} 
(2005), no. 3, 471--490.

\noindent [IkKeMu2] I. Ikromov, M. Kempe, and D. M\"uller, {\it Estimates for maximal functions associated
to hypersurfaces in $R^3$ and related problems of harmonic analysis}, Acta Math. {\bf 204} (2010), no. 2,
151--271.

\noindent [IkMu] I. Ikromov, D. M\"uller, {\it On adapted coordinate systems}, to appear, Trans. AMS.

\noindent [IoSa1] A. Iosevich, E. Sawyer, {\it Oscillatory integrals and maximal averages over homogeneous
surfaces}, Duke Math. J. {\bf 82} no. 1 (1996), 103-141.

\noindent [IoSa2] A. Iosevich, E. Sawyer, {\it Maximal averages over surfaces},  Adv. Math. {\bf 132} 
(1997), no. 1, 46--119.

\noindent [NaSeWa] A. Nagel, A. Seeger, and S. Wainger, {\it Averages over convex hypersurfaces},
Amer. J. Math. {\bf 115} (1993), no. 4, 903--927.

\noindent [PSt] D. H. Phong, E. M. Stein, {\it The Newton polyhedron and oscillatory integral operators}, 
Acta Math. {\bf 179} (1997), 107-152.

\noindent [So] C. Sogge, {\it Maximal operators associated to hypersurfaces with one nonvanishing principal 
curvature}  (English summary)  in {\it Fourier analysis and partial differential equations} (Miraflores de 
la Sierra, 1992),  317--323, Stud. Adv. Math., CRC, Boca Raton, FL, 1995. 

\noindent [SoSt] C. Sogge and E. Stein, {\it Averages of functions over hypersurfaces in $R^n$}, Invent.
Math. {\bf 82} (1985), no. 3, 543--556.

\noindent [St1] E. Stein, {\it Maximal functions. I. Spherical means.} Proc. Nat. Acad. Sci. U.S.A. 
{\bf 73} (1976), no. 7, 2174--2175. 

\noindent [St2] E. Stein, {\it Harmonic analysis; real-variable methods, orthogonality, and oscillatory 
integrals}, Princeton Mathematics Series Vol. 43, Princeton University Press, Princeton, NJ, 1993.

\noindent [V] A. N. Varchenko, {\it Newton polyhedra and estimates of oscillatory integrals}, Functional 
Anal. Appl. {\bf 18} (1976), no. 3, 175-196.

\line{}
\line{}

\noindent Department of Mathematics, Statistics, and Computer Science \hfill \break
\noindent University of Illinois at Chicago \hfill \break
\noindent 322 Science and Engineering Offices \hfill \break
\noindent 851 S. Morgan Street \hfill \break
\noindent Chicago, IL 60607-7045 \hfill \break
\end